\documentclass[11pt,a4paper]{article}

\usepackage{amsmath,amsfonts,amssymb,amsthm}

\newcommand{\diam}{\operatorname{diam}}

\newcommand{\counte}{equation}%%%设置定理与公式公用计数器
\newtheorem{theorem}[\counte]{\bf Theorem}%[section]

\newtheorem{defn}[\counte]{\bf Definition}

\newtheorem{prop}[\counte]{\bf Proposition}
\newtheorem{lemma}[\counte]{\bf Lemma}

\newtheorem{pf}[\counte]{\bf Proof}
\newtheorem{coro}[\counte]{\bf Corollary}

\newtheorem{remark}[\counte]{\bf Remark}

\allowdisplaybreaks

\numberwithin{equation}{section}%%公式随section重新计数
\renewcommand{\thefootnote}{\fnsymbol{footnote}}

\begin{document}

\renewcommand{\thefootnote}{\arabic{footnote}}

\centerline{\bf On $\frac\pi2$-separated subsets of Alexandrov spaces with curvature $\geq1$
\footnote{Supported by NSFC 11001015 and 11171025.
\hfill{$\,$}}}

\vskip5mm

\centerline{\it Dedicated to Xiaochun Rong for his 60th birthday}

\vskip5mm

\centerline{ Xiaole Su, Hongwei Sun, Yusheng Wang\footnote{The corresponding author (E-mail: wwyusheng@gmail.com). \hfill{$\,$}}}

\vskip6mm

\noindent{\bf Abstract.} Let $M$ be an $n$-dimensional Alexandrov space with curvature $\geq 1$,
and let $\{q_1,\cdots,q_k\}$ be any $\frac\pi2$-separated subset in $M$
(i.e. the distance $|q_iq_j|\geq\frac{\pi}{2}$ for any $i\neq j$).
Under the additional conditions ``$|q_iq_j|<\pi$'' and ``the diameter $\diam(M)\leq \frac\pi2$'',
we respectively give the upper bound of $k$ (which depends only on $n$),
and we classify the (topological or geometric) structure
of $M$ when $k$ attains the upper bound.

\vskip1mm


\vskip1mm


\vskip6mm

\setcounter{section}{-1}

\section{Introduction}

In studying the Morse theory on
Alexandrov spaces with curvature $\geq \kappa$ ([P1]), the following basic and easy idea plays an important role.

\begin{theorem}\label{thm0.1}   {\it
Let $\{q_1,\cdots,q_k\}$ be a subset in an $n$-dimensional
Alexandrov space with curvature $\geq 1$. If the distance
$|q_iq_j|>\frac{\pi}{2}$ for any $1\leq i\neq j\leq k$, then $k\leq
n+2$}.
\end{theorem}

\noindent However, Theorem \ref{thm0.1} is not formulated in [P1]
(maybe due to the simplicity of it, especially to Perel$'$man). On the $M$ in Theorem \ref{thm0.1},
in [GW] there are the following two packing radius theorems when $k\geq 2$: one is that {\it $M$ is homeomorphic to
the join $\Bbb S^{k-2}*N$};
the other is that {\it $\min_{i\neq j}\{|q_iq_j|\}\leq\arccos(\frac{-1}{k-1}),$ and if the equality holds
then $M$ is isometric to the join $\Bbb S^{k-2}*N$, where $\Bbb S^{k-2}$ is the $(k-2)$-dimensional unit sphere
and $N$ is some $(n-k+1)$-dimensional Alexandrov space with curvature $\geq 1$}. (See the comments after Theorem B below
for the definition of the join.)

We find that if ``$|q_iq_j|>\frac{\pi}{2}$'' is changed to
``$|q_iq_j|\geq\frac{\pi}{2}$'' in Theorem \ref{thm0.1}, to find the upper bounds of $k$
under the condition $$|q_iq_j|<\pi$$ is more interesting and
difficult. In the present paper, we make clear this upper bound, and we can classify the geometric structure of $M$
if $k$ attains the upper bound. Note that the condition ``$|q_iq_j|>\frac{\pi}{2}$'' in Theorem \ref{thm0.1} implies that
$|q_iq_j|<\pi$ if $k\geq 3$ (see Lemma \ref{lem1.2} below), which is an important idea in [P1]. Hence, the condition ``$|q_iq_j|<\pi$''
is not an artificial one.

Of course, if we only change ``$|q_iq_j|>\frac{\pi}{2}$'' to
``$|q_iq_j|\geq\frac{\pi}{2}$'' in Theorem \ref{thm0.1},
we have the following well-known result (cf. [GW]).

\begin{theorem}\label{thm0.2}   Let $\{q_1,\cdots,q_k\}$ be a subset in an $n$-dimensional
Alexandrov space with curvature $\geq 1$. If the distance
$|q_iq_j|\geq\frac{\pi}{2}$ for any $1\leq i\neq j\leq k$, then
$k\leq2(n+1)$, and if the equality holds then $M$ is isometric to
$\Bbb S^n$ and we can rearrange all $q_i$ such that
$|q_{2j-1}q_{2j}|=\pi$ for $1\leq j\leq n+1$.
\end{theorem}

\noindent For the completeness of the paper, we will give a proof of
Theorem \ref{thm0.2} in Appendix.

\vskip2mm

In the present paper, we let ${\cal A}^n(\kappa)$  denote the
collection of all $n$-dimensional Alexandrov spaces with curvature
$\geq \kappa$ (containing all $n$-dimensional Riemannian manifolds
with sectional curvature $\geq \kappa$), and without special remark
we always consider complete spaces in ${\cal A}^n(\kappa)$.

\begin{defn}{\rm  Let $M\in {\cal A}^n(1)$, and let
$Q\triangleq\{q_1,\cdots,q_k\}$ be a subset in $M$. We call $Q$ a {\it $\frac\pi2$-separated subset} in $M$ if
the distance $|q_iq_j|\geq\frac{\pi}{2}$ for any $1\leq i\neq j\leq k$.}
\end{defn}

Now we give our first estimate result.

\vskip2mm

\noindent{\bf Theorem A}  {\it Let $M\in {\cal A}^n(1)$, and let
$\{q_1,\cdots,q_k\}$ be a $\frac\pi2$-separated subset in $M$. If $|q_1q_i|>\frac\pi2$ for any
$2\leq i\leq k$, then $k\leq n+2$; and if the equality holds, then $M$
is homeomorphic to $\Bbb S^n$ (and thus $M$ has empty boundary).}

\vskip2mm

Note that Theorem A implies Theorem 0.1. Since the idea of
estimating $k$ in Theorem A is the same as in Theorem 0.1, the upper
bound of $k$ in Theorem A should be known to experts. For the
convenience of readers, we will give its proof in Section 1.
However, the following results are not so obvious.

\vskip2mm

\noindent{\bf Theorem B} {\it  Let $M\in {\cal A}^n(1)$, and let
$\{q_1,\cdots,q_k\}$ be a $\frac\pi2$-separated subset in $M$. If
$|q_iq_j|<\pi$ for any $1\leq i\neq j\leq k$, then  $$k\leq 3l\
(\text{resp. } 3l+1) \text{ for } n=2l-1\ (\text{resp. } 2l);$$
moreover, if the equality holds, then we can rearrange all $q_i$ such
that $M$ is isometric to $S_1^1*\cdots *S_l^1$ (resp. if $M$ has
empty boundary, then either $M$ is isometric to $S_1^1*\cdots
*S_{l-1}^1*N$ for some $N\in {\cal A}^2(1)$,  or
$\{q_{3l+1}\}*S_1^1*\cdots *S_l^1$ can be isometrically embedded
into $M$; if $M$ has nonempty boundary, then $M$ is isometric to
$\{q_{3l+1}\}*S_1^1*\cdots *S_l^1$) with $S_j^1$ having perimeter
$\geq\frac{3\pi}{2}$ (of course $\leq 2\pi$) and $q_{3j-2}, q_{3j-1}, q_{3j}\in S_j^1$ for
each $j$.}

\vskip2mm

Recall that a join $X*Y$ with $X, Y\in{\cal A}(1)$ is defined as
follows ([BGP]). $X*Y=X\times Y\times[0,\frac{\pi}{2}]/\sim$, where
$(x,y,t)\sim(x',y',t')\Leftrightarrow t=t'=0 \text{ and } x=x' \text{ or }
t=t'=\frac{\pi}{2}\text{ and } y=y'$,
and for any $p_i=[(x_i,y_i,t_i)]\in X*Y$
$$\cos|p_1p_2|=\cos t_1\cos t_2\cos|x_1x_2|+\sin t_1\sin t_2\cos|y_1y_2|.$$
Note that $X*Y$ also belongs to ${\cal A}(1)$ and
$\dim(X*Y)=\dim(X)+\dim(Y)+1$, and $X*Y$ is a Riemannian manifold if and
only if $X$ and $Y$ are isometric to unit spheres.

A very interesting corollary of Theorem B is on Riemannian cases.

\vskip2mm

\noindent{\bf Corollary C} {\it   Let $M$ be a closed $n$-dimensional
Riemannian manifold with sectional curvature $\geq 1$, and let
$\{q_1,\cdots,q_k\}$ be a $\frac\pi2$-separated subset in $M$ with
$|q_iq_j|<\pi$ for any $1\leq i\neq j\leq k$. Then $k\leq 3l\
(\text{resp. } 3l+1) \text{ for } n=2l-1\ (\text{resp. } 2l)$;
and if the equality holds and $n>2$, then $M$ is isometric to the unit sphere $\Bbb S^n$.}

\vskip2mm

When $n=3$, for example, Corollary C says that $M$ contains at most
6 points $q_1,\cdots,q_6$ with $\frac\pi2\leq |q_iq_j|<\pi$ for any
$1\leq i\neq j\leq 6$, and only the unit sphere $\Bbb S^3$ contains
such 6 points (if we embed $\Bbb S^3$ isometrically into  the
Euclidean space $\Bbb R^4=\{(x_1,x_2,x_3,x_4)|x_i\in \Bbb R\}$, we
can select the former (resp. latter) 3 points on the plane
$\{(x_1,x_2,0,0)\}$ (resp. $\{(0,0,x_3,x_4)\}$). And this is the
unique way to select such 6 points up to an orthogonal
transformation of $\Bbb R^4$.).

\vskip2mm

\noindent{\bf Theorem D}  {\it Let $M\in {\cal A}^n(1)$, and let
$\{q_1,\cdots,q_k\}$ be a $\frac\pi2$-separated subset in $M$. If the
diameter $\diam(M)\leq \frac{\pi}{2}$, then there exists an
isometrical embedding $f: \Delta^{k-1}_+\to M$ such that
$q_1,\cdots,q_k$ are the vertices of $f(\Delta^{k-1}_+)$, where
$$\Delta^{k-1}_+\triangleq\left\{(x_1,\cdots, x_k)\in \mathbb{R}^{k}|
\sum\limits_{i=1}^{k} x_i^2=1, x_i\geq 0\right\}\subset
\mathbb{S}^{k-1}.$$ As a result, $k\leq n+1$; moreover, if $k=n+1$,
then $M$ is a glued space of finite copies of $\Delta^{n}_+$ along
some ``faces'' $\Delta^{n-1}_+$ of them.}

\vskip2mm

We know that the boundary of a $\Delta^{n}_+$ consists of $n+1$
copies of $\Delta^{n-1}_+$. Here, such a $\Delta^{n-1}_+$ is said
to be a ``{\it face}'' of the $\Delta^{n}_+$.

Similarly, Theorem D has the following corollary on Riemannian manifolds.

\vskip2mm

\noindent{\bf Corollary E} {\it If in addition $M$ is a closed $n$-dimensional
Riemannian manifold with sectional curvature $\geq 1$ in Theorem D, and if $k=n+1$,
then $M$ is isometric to the projective space $\Bbb R\Bbb P^n$ with the canonical metric
(i.e. the metric induced from $\Bbb S^n$).}

\vskip2mm

On Theorem D, we supply another Riemannian example (for more general examples please
refer to Remark \ref{rem3.12}). We consider the
complex projective space $\Bbb C\Bbb P^n$ with the canonical metric
(i.e. the metric induced from $\Bbb S^{2n+1}$). It is well known that
$\Bbb C\Bbb P^n$ has sectional curvature $\geq 1$ (and $\leq4$) and the
diameter $\leq\frac\pi2$.  By the induction on $n$,
it is not hard to see that $\Bbb C\Bbb P^n$ contains $\{q_1,\cdots,q_{n+1}\}$ with
$|q_iq_j|=\frac\pi2$ for any $1\leq i\neq j\leq n+1$
(however, we cannot find a $\frac\pi2$-separated subset containing $n+2$ points in $\Bbb C\Bbb P^n$). According to Theorem D,
$\Delta^{n}_+$ can be isometrically embedded into $\Bbb C\Bbb P^n$.

\vskip2mm

We will end this section by introducing our mail tool---the Toponogov Comparison
Theorem, which is the essential geometry in Alexandrov spaces with curvature $\geq\kappa$.

We always let $[pq]$ denote a geodesic (i.e. a shortest path) between $p$
and $q$ in $M\in\mathcal{A}^n(\kappa)$, and let $\uparrow_p^q$ denote the direction at $p$ of the
geodesic $[pq]$. Given another geodesic $[pr]$ in $M$, we let $\angle qpr$ denote the angle between
$[pq]$ and $[pr]$ at $p$ (for the detailed contents of angles please refer to [BGP]). We know that
$\angle qpr$ is equal to $|\uparrow_p^q\uparrow_p^r|$
(i.e. the distance between  $\uparrow_p^q$ and $\uparrow_p^r$) in $\Sigma_pM$, where $\Sigma_pM \in \mathcal{A}^{n-1}(1)$ is the direction space of $M$ at $p$.

We say that the geodesics $[pq]$ and $[pr]$  form a {\it hinge}
$p\prec^q_r$ at $p$ with angle $\angle qpr$, and call an
associated hinge $\tilde p\prec^{\tilde q}_{\tilde r}$ in $\Bbb
S^2_\kappa$ with  $|\tilde p\tilde q|=|pq|$, $|\tilde p\tilde
r|=|pr|$ and $\angle\tilde q\tilde p\tilde r=\angle qpr$ the
comparison hinge of $p\prec^q_r$, where $\Bbb S^2_\kappa$ is the
complete and simply-connected 2-manifold of constant curvature
$\kappa$. Similarly, we say that geodesics $[pq], [qr]$ and $[rp]$
form a {\it triangle} $\triangle pqr$, and call an associated
triangle $\triangle\tilde p\tilde q\tilde r$ in $\Bbb S^2_\kappa$
with $|\tilde p\tilde q|=|pq|,|\tilde p\tilde r|=|pr|$ and $|\tilde
r\tilde q|=|rq|$ the comparison triangle of $\triangle pqr$.

For any triangle $\triangle pqr$ (we only need to consider the case
$|pq|+|pr|+|qr|<2\pi/\sqrt\kappa$ if $\kappa>0$ ([BGP])) and hinge
$p\prec^q_r$ in $M\in\mathcal{A}^n(\kappa)$ and their comparison
triangle and hinge $\triangle\tilde p\tilde q\tilde r$ and $\tilde
p\prec^{\tilde q}_{\tilde r}$, the Toponogov Comparison Theorem
(TCT) asserts that ([BGP]):

\begin{theorem}[TCT]\label{tct}

\noindent {\rm(i)} For any two points $s\in[qr]\subset \triangle pqr$ and
$\tilde s\in[\tilde q\tilde r]\subset\triangle\tilde p\tilde q\tilde
r$ with $|qs|=|\tilde q\tilde s|$, we have $|ps|\geq|\tilde p\tilde
s|$.

\vskip1mm

\noindent {\rm(ii)}  In $\triangle pqr$ and $\triangle\tilde p\tilde
q\tilde r$, we have  $\angle pqr\geq\angle \tilde p\tilde q\tilde
r,$ $\angle qrp\geq\angle \tilde q\tilde r\tilde p$ and $\angle
rpq\geq\angle \tilde r\tilde p\tilde q$.

\vskip1mm

\noindent {\rm(iii)}  In $p\prec^q_r$ and $\tilde
p\prec^{\tilde q}_{\tilde r}$, we have  $|\tilde q\tilde
r|\geq|qr|$.

\end{theorem}

It is known that (i)-(iii) of TCT are equivalent to each other.
Moreover, we have the following result when the ``$=$'' holds in
TCT.

\begin{theorem}[TCT for ``$=$'' ({[GM]})] \label{tct=}  {\rm(i)} If there is a point
$s\in[qr]^\circ$ such that $|ps|=|\tilde p\tilde s|$ in {\rm (i)} of
TCT, then for any given geodesic $[ps]$ there exist unique two
geodesics $[pq]'$ and $[pr]'$ (maybe not $[pq]$ and $[pr]$) such
that the triangle formed by $[pq]'$, $[pr]'$ and $[qr]$ is isometric
to its comparison triangle.

\noindent{\rm(ii)} If $|\tilde q\tilde r|=|qr|$ in {\rm (iii)} of
TCT (or if $\angle rpq=\angle \tilde r\tilde p\tilde q$ in {\rm
(ii)} of TCT), then there exists geodesic $[qr]'$ (maybe not $[qr]$)
such that the triangle formed by $[pq]$, $[pr]$ and $[qr]'$ is
isometric to its comparison triangle.
\end{theorem}

\vskip2mm

The rest of the paper is organized as follows. In Sections 1-3, we give the proofs of Theorem A,
Theorem B and Corollary C, and Theorem D and Corollary E respectively. A technical corollary of Theorem D
is given in Section 4. In Appendix, we will prove Theorem \ref{thm0.2} and Lemmas \ref{lem1.4} and \ref{lem2.11}.

%%%%%%%%%%%%%%%%%%%%%%%%%%%%%%%%%%%%%%% Section 1  %%%%%%%%%%%%%%%%%%%%%%%%%%%%%%%%%%%%%%%

\section{Proof of Theorem A}

In the paper, we often use the following lemma,
an obvious corollary of Theorem \ref{tct}.

\begin{lemma}\label{lem1.1}
Let $M\in {\cal A}^n(1)$,  and let $p, q, r\in M$ with
$|qr|\geq\frac{\pi}{2}$. If either $|pq|,|pr|\leq\frac\pi2$ or
$\frac\pi2\leq|pq|,|pr|<\pi$, then for any geodesics $[pq]$ and
$[pr]$ we have $|\uparrow_{p}^{q}\uparrow_{p}^{r}|\geq\frac\pi2$ in
$\Sigma_{p}M$; and if in addition $|qr|>\frac\pi2$ or $|pq|,|pr|>\frac\pi2$, then
$|\uparrow_{p}^{q}\uparrow_{p}^{r}|>\frac\pi2$.
\end{lemma}

And the following basic fact will be used sometimes.

\begin{lemma}[{[BGP]}]\label{lem1.2}
Let $M\in {\cal A}^{n}(1)$ and $p,q\in M$. If $|pq|=\pi$, then $|px|+|qx|=\pi$ for any $x\in M$, and $M=\{p,q\}*M_1$ for some $M_1\in {\cal A}^{n-1}(1)$.
\end{lemma}

Let $M\in {\cal A}^n(1)$. For $n=0$ and 1, we make the following convention: if $n=0$, then $M$ consists of one point or two points with distance equal to $\pi$;
if $n=1$, then $M$ is an arc with length $\leq \pi$ or a circle with perimeter $\leq 2\pi$.

\vskip2mm

\noindent{\it Proof of Theorem A}.

We will give the proof by the induction on the dimension $n$.
Obviously, Theorem A is true if $n=0$ and 1 (see the above
convention). Now we assume that $n>1$, and we can assume that
$k\geq3$. According to Lemma \ref{lem1.2}, ``$|q_1q_i|>\frac\pi2$
for any $2\leq i\leq k$'' implies that $|q_iq_j|<\pi$ for any $1\leq
i\neq j\leq k$. Then by Lemma \ref{lem1.1}, any
$\{\uparrow_{q_k}^{q_1},\cdots,\uparrow_{q_k}^{q_{k-1}}\}$ is a
$\frac\pi2$-separated subset in $\Sigma_{q_k}M\in {\cal A}^{n-1}(1)$
with $|\uparrow_{q_k}^{q_1}\uparrow_{q_k}^{q_{i}}|>\frac\pi2$ for
all $2\le i\le k-1$. By the inductive assumption on $\Sigma_{q_k}M$,
we have $$k-1\leq n-1+2, \text{ i.e., } k\leq n+2.$$

Now we will prove that $M$ is homeomorphic to $\Bbb S^n$ if $k=n+2$.
By the Radius Sphere Theorem ([GP]), it suffices to show that
$\text{rad}(M)>\frac\pi2$, where $\text{rad}(M)$ is the radius of
$M$ defined by $\min_{p\in M}\{\max_{q\in M}|pq|\}$. Note that if
$\text{rad}(M)\le\frac\pi2$, then there is a point $p\in M$ such
that $|px|\leq\frac{\pi}{2}$ for all $x\in M$. Obviously,
$p\not\in\{q_1,\cdots,q_{k}\}$ (note that $|q_1q_i|>\frac\pi2$ for
all $2\le i\le k$). Then by Lemma \ref{lem1.1}, any
$\{\uparrow_{p}^{q_1},\cdots,\uparrow_{p}^{q_{k}}\}$ is a
$\frac\pi2$-separated subset in $\Sigma_{p}M\in {\cal A}^{n-1}(1)$
with $|\uparrow_{p}^{q_1}\uparrow_{p}^{q_{i}}|>\frac\pi2$ for all
$2\le i\le k$, so by the former part (we have proved) we have $k\leq
n-1+2=n+1$; a contradiction. It therefore has to hold that the
radius $\text{rad}(M)>\frac\pi2$ (and thus $M$ is homeomorphic to
$\Bbb S^n$). \hfill$\Box$

\vskip2mm

In the above proof, we use the Radius Sphere Theorem to show that
$M$ is homeomorphic to $\Bbb S^n$, and thus $M$ has empty boundary.
In fact, we can prove that $M$ has empty boundary (when $k=n+2$ in
Theorem A) without the Radius Sphere Theorem as follows.

\vskip2mm

\begin{pf}[a proof for ``$M$ has empty boundary if $k=n+2$ in Theorem A'']\label{pf1.3}

\rm \hskip2cm

Obviously, this is true when $n=1$. Next, we will
derive a contradiction by applying the induction on $n$ and assuming that
the boundary $\partial M\neq \emptyset$.

We consider $\Sigma_{q_2}M$ ($\in{\cal A}^{n-1}(1)$).
From the above proof, we know that $|q_iq_j|<\pi$ for any $1\leq i\neq j\leq n+2$.
By Lemma \ref{lem1.1}, any
$\{\uparrow_{q_2}^{q_1},\uparrow_{q_2}^{q_3},\cdots,\uparrow_{q_2}^{q_{n+2}}\}$
is a $\frac\pi2$-separated subset with
$|\uparrow_{q_2}^{q_1}\uparrow_{q_2}^{q_i}|>\frac{\pi}{2}$ for any
$3\leq i\leq n+2$. By the inductive assumption, $\Sigma_{q_2}M$ has empty
boundary, so $q_2\not\in \partial M$ (for the detailed contents
on the boundary of a space in ${\cal A}^{n}(\kappa)$ please refer to [BGP]).
Then we select $p\in\partial M$ such that $|q_2p|=|q_2\partial M|$. If
$|q_2p|\geq\frac\pi2$, then by Lemma \ref{lem1.4} below
$M=\{q_2\}*\partial M$, which contradicts
``$|q_2q_1|>\frac{\pi}{2}$''. Now we can assume that
$|q_2p|<\frac\pi2$. On the other hand, since $|q_2p|\leq|q_2x|$ for
all $x\in\partial M$, by the first variation formula ([BGP]) we have
$$|\uparrow_p^{q_2}\xi|\geq\frac\pi2$$
for any geodesic $[pq_2]$ and $\xi\in
\partial(\Sigma_pM)$ (refer to [BGP] for $\partial(\Sigma_pM)$).
By Lemma \ref{lem1.4} below,
$$\Sigma_pM=\{\uparrow_p^{q_{2}}\}*\partial(\Sigma_pM).$$
Hence, for any geodesic $[pq_i]$ with $i\neq 2$, we have
$|\uparrow_p^{q_{2}}\uparrow_p^{q_i}|\leq\frac\pi2$. Due to Theorem \ref{tct},
we can conclude that $$|pq_i|\geq\frac\pi2\ (i\neq 2)$$ by considering
the comparison triangle of $\triangle q_2pq_i$ containing sides
$[pq_2]$ and $[pq_i]$ (note that $|pq_2|<\frac\pi2$, $|q_2q_i|\geq\frac\pi2$ and $\angle q_2pq_i\leq\frac\pi2$).
Furthermore, by Lemma \ref{lem1.2} again, we can conclude that $|pq_i|<\pi$ for $i\neq 2$
(because $|q_1q_j|>\frac\pi2$ for $j=3,\cdots, n+2$). Hence, by Lemma \ref{lem1.1}, any
$\{\uparrow_{p}^{q_1},\uparrow_{p}^{q_3},\cdots,\uparrow_{p}^{q_{n+2}}\}$
is a $\frac\pi2$-separated subset in $\Sigma_pM\in {\cal A}^{n-1}(1)$ with
$|\uparrow_{p}^{q_1}\uparrow_{p}^{q_j}|>\frac{\pi}{2}$ for any
$3\leq j\leq n+2$. By the inductive assumption, $\Sigma_{p}M$ has empty
boundary, which contradicts ``$p\in \partial M$''.\hfill$\Box$
\end{pf}

\begin{lemma}\label{lem1.4}
Let $M\in {\cal A}^n(1)$ with nonempty boundary. If $|p\partial
M|\geq\frac\pi2$ for some $p\in M$, then $M=\{p\}*\partial M$.
\end{lemma}

It is easy to see that Lemma \ref{lem1.4} is a corollary of the
Doubling Theorem by Perel$'$man ([P2]). For the convenience of readers, we
will give an elementary proof for it in Appendix.

%%%%%%%%%%%%%%%%%%%%%%%%%%%%%%%%%%%%%%% Section 2  %%%%%%%%%%%%%%%%%%%%%%%%%%%%%%%%%%%%%%%

\section{Proofs of Theorem B and Corollary C}

We will prove the following generalized version of Theorem B.

\begin{theorem} \label{thm2.1}
Let $M\in {\cal A}^n(1)$, and let
$\{q_1,\cdots,q_h,q_{h+1},\cdots,q_k\}$ be a $\frac\pi2$-separated
subset in $M$. Suppose that $|q_iq_j|<\pi$ for any $1\leq i\neq
j\leq k$, and that $|q_1q_i|>\frac\pi2$ for any $2\leq i\leq
h$. Then $h\leq n+2$, and the following hold:

\noindent {\rm (i)} If $h=n+2$, then $k=n+2$, and $M$ has empty
boundary.

\noindent {\rm (ii)} If $h=n+1$, then $k\leq n+2$; and if the equality holds,
then either $M$ has empty boundary, or $M=\{q_{n+2}\}*N$ for some
$N\in {\cal A}^{n-1}(1)$ without boundary.

\noindent {\rm (iii)} If $4\leq h\leq n$, then $k-h\leq 3l$ (resp.
$3l+1$) for $n-h+1=2l-1$ (resp. $2l$); and if the equality holds, then
$M$ is isometric to $L*S_1^1*\cdots*S_l^1$ (resp. either $M$ is isometric
to $N*S_1^1*\cdots* S_{l-1}^1$, or $L*S_1^1*\cdots*S_l^1*\{q_i\}$
for some $i>h$ can be isometrically embedded into $M$), where $S^1_j$ is
of perimeter $\geq\frac{3\pi}{2}$, $L\in {\cal A}^{h-2}(1)$ and
$N\in {\cal A}^{h+1}(1)$.

\noindent {\rm (iv)}  If $h\leq 3$, then $k\leq 3l$ (resp. $3l+1$)
for $n=2l-1$ (resp. $2l$); moreover, if the equality holds, then  we can
rearrange all $q_i$ such that $M$ is isometric to $S_1^1*\cdots
*S_l^1$ (resp. if $M$ has empty boundary, then either $M$ is isometric to
$S_1^1*\cdots *S_{l-1}^1*N$ for some $N\in {\cal A}^2(1)$,  or
$\{q_{3l+1}\}*S_1^1*\cdots *S_l^1$ can be isometrically embedded
into $M$; if $M$ has nonempty boundary, then $M$ is isometric to
$\{q_{3l+1}\}*S_1^1*\cdots *S_l^1$) with $S_j^1$ having perimeter
$\geq\frac{3\pi}{2}$ (of course $\leq2\pi$) and $q_{3j-2}, q_{3j-1}, q_{3j}\in S_j^1$ for
each $j$.
\end{theorem}

Obviously, the conclusion  ``$h\leq n+2$'' in Theorem \ref{thm2.1}
is included in Theorem A. And note that Theorem B is
included in (iv) of Theorem \ref{thm2.1}.

In the following we will first give the proofs of (i) and (ii) in Theorem \ref{thm2.1}.

\vskip1mm

\noindent{\bf Proof of (i) in Theorem \ref{thm2.1}}:

By Theorem A, $M$ has empty boundary, so we only need to show that
$k=n+2$. If $k>n+2$, then we consider $\Sigma_{q_{h+1}}M\in{\cal
A}^{n-1}(1)$. By Lemma \ref{lem1.1}, any
$\{\uparrow_{q_{h+1}}^{q_1},$ $\cdots, \uparrow_{q_{h+1}}^{q_{h}}\}$ is a
$\frac\pi2$-separated subset in $\Sigma_{q_{h+1}}M$ with
$|\uparrow_{q_{h+1}}^{q_1}\uparrow_{q_{h+1}}^{q_{i}}|>\frac\pi2$ for
any $2\leq i\leq h$. By Theorem A, we have $h\leq n-1+2$, which
contradicts $h=n+2$.

\hfill$\Box$

\vskip1mm

\noindent{\bf Proof of (ii) in Theorem \ref{thm2.1}}:

Obviously, this is true if $n=0$ and 1. Then we assume that $n\geq 2$,
which implies that $h\geq 3$.

We first prove that $k\leq n+2$. If $k>n+2$, then in
$\Sigma_{q_{h+2}}M\in{\cal A}^{n-1}(1)$ any
$\{\uparrow_{q_{h+2}}^{q_1},$ $\cdots,\uparrow_{q_{h+2}}^{q_{h}},$
$\uparrow_{q_{h+2}}^{q_{h+1}}\}$ is a $\frac\pi2$-separated subset
with $|\uparrow_{q_{h+2}}^{q_1}\uparrow_{q_{h+2}}^{q_{i}}|>\frac\pi2$ for
any $2\leq i\leq h$ (by Lemma
\ref{lem1.1}). Furthermore, by Lemma \ref{lem1.2} we can conclude that
$|\uparrow_{q_{h+2}}^{q_i}\uparrow_{q_{h+2}}^{q_{j}}|<\pi$ for any
$1\leq i\neq j \leq h+1$.
Therefore, the $\frac\pi2$-separated subset
$\{\uparrow_{q_{h+2}}^{q_1},\cdots,\uparrow_{q_{h+2}}^{q_{h}},\uparrow_{q_{h+2}}^{q_{h+1}}\}$
in $\Sigma_{q_{h+2}}M\in{\cal A}^{n-1}(1)$ satisfies the conditions
of Theorem \ref{thm2.1}. Then by (i) of Theorem
\ref{thm2.1}, it has to hold that $h<(n-1)+2=n+1$ which contradicts $h=n+1$.

\vskip1mm

Next we only need to prove that if $k=n+2$ and if $M$ has nonempty
boundary, then $M=\{q_{n+2}\}*\partial M$. By Lemma \ref{lem1.4}, it
suffices to show that $|q_{n+2}\partial M|\geq\frac\pi2$. If
$|q_{n+2}\partial M|<\frac\pi2$, we select $p\in
\partial M$ such that $|q_{n+2}p|=|q_{n+2}\partial M|$. Then like Proof \ref{pf1.3},
we can get that $\frac\pi2\leq|pq_i|<\pi$
for $1\leq i\leq n+1$. Hence, by Lemma \ref{lem1.1} any
$\{\uparrow_{p}^{q_1},\cdots,\uparrow_{p}^{q_{n+1}}\}$
is a $\frac\pi2$-separated subset in $\Sigma_pM\in {\cal A}^{n-1}(1)$ with
$|\uparrow_{p}^{q_1}\uparrow_{p}^{q_i}|>\frac{\pi}{2}$
for any $2\leq i\leq n+1$. By (i) of Theorem \ref{thm2.1},
$\Sigma_{p}M$ has empty boundary, which contradicts $p\in \partial
M$. \hfill$\Box$

\subsection{Some preparations for proving (iii) and (iv) of Theorem \ref{thm2.1}}

Given a subset $A$ of $M$, we let $A^{\geq
d}\triangleq\{x\in M||xa|\geq d, \forall\ a\in A\}$. And similarly we can define the corresponding $A^{\leq
d}$, $A^{=d}$, $A^{<d}$ and $A^{>d}$.
From Theorem \ref{tct}, we can immediately see the following lemma.

\begin{lemma}\label{lem2.2}
Let $A$ be a subset of $M\in {\cal A}^n(1)$. Then $A^{\geq
\frac{\pi}{2}}$ is  convex in $M$.
\end{lemma}

Recall that $N$ is said to be {\it convex} in $M\in {\cal A}^n(1)$
if there is a geodesic $[xy]$ belonging to $N$ for any $x,y\in N$,
or $N$ consists of two points with distance equal to $\pi$, or $N$
consists of only one point. We know that a convex subset $N$ in $M$
also belongs to ${\cal A}^m(1)$; and if $N\subsetneq M$ and $N$ has
empty boundary, then $m<n$ ([BGP]).

\begin{lemma}\label{lem2.3} {\bf ([Y])}
Let $M\in {\cal A}^n(1)$, and let $A$ be a complete locally convex subset in
$M$. If $A$ has empty boundary, then $A^{\geq
\frac{\pi}{2}}=A^{=\frac{\pi}{2}}$.
\end{lemma}

In our proof, we will use a special and generalized case of Lemma \ref{lem2.3}.

\begin{lemma}\label{lem2.4}
Let $M\in {\cal A}^n(1)$, and let
$A\triangleq\bigcup_{i=0}^l[p_ip_{i+1}]\subset M$ (where
$p_{l+1}=p_0$). If geodesics $\{[p_ip_{i+1}]\}_{i=0}^l$ satisfy
$|\uparrow_{p_i}^{p_{i-1}}\uparrow_{p_i}^{p_{i+1}}|=\pi$ (where
$p_{-1}=p_l$), then $A^{\geq \frac{\pi}{2}}=A^{=\frac{\pi}{2}}$.
\end{lemma}

From Lemma \ref{lem2.5} below, we have that
$\dim(A)+\dim(A^{=\frac{\pi}{2}})\leq n-1$ in Lemmas \ref{lem2.3} and \ref{lem2.4}.

\begin{lemma}\label{lem2.5} {\bf ([RW])}
Let $M\in {\cal A}^n(1)$, and let $N_1$ and $N_2$ be two locally convex
subsets in $M$. If $|x_1x_2|=\frac\pi2$ for all $x_i\in N_i$, then
$\dim(N_1)+\dim(N_2)\leq n-1$.
\end{lemma}

Based on Lemmas \ref{lem2.2}-\ref{lem2.5}, we will give the proofs
of (iii) and (iv) in Theorem \ref{thm2.1} by the induction on $n$.
Because the inductive processes of the proofs for (iii) and (iv) are
almost identical (please see Remark \ref{rem2.10} below for the main
difference between them), we only give the detailed proof for (iv).
Obviously, (iv) is true when $n=1$.

\subsection{Estimating $k$ for $n=2$ and $3$ in (iv) of Theorem \ref{thm2.1}}

In this subsection, we will mainly prove that $k\leq 4$ and $6$ when $n=2$ and $3$ respectively in (iv) of Theorem \ref{thm2.1}. In the proof, we need the following technique lemma.

\begin{lemma}\label{lem2.6}  Let $\{q_1,q_2,q_3\}$ be a
$\frac\pi2$-separated subset in $M\in {\cal A}^n(1)$. If
$|q_1q_i|<\pi$ and there are geodesics $[q_1q_i]$ ($i=2, 3$)
such that $|\uparrow_{q_1}^{q_2}\uparrow_{q_1}^{q_3}|=\pi$, then
$\{q_1,q_2,q_3\}^{\geq\frac\pi2}\subseteq\{q_2,q_3\}^{=\frac\pi2}$;
and for any $z\in
\{q_1,q_2,q_3\}^{\geq\frac\pi2}\cap\{q_1\}^{<\pi}$, there is a
geodesic $[q_iz]$ such that
$\frac\pi2\leq|\uparrow_{q_i}^{z}\uparrow_{q_i}^{q_1}|<\pi$ ($i=2, 3$).
\end{lemma}

\noindent{\it Proof}. We first note that if $z\in
\{q_1,q_2,q_3\}^{\geq\frac\pi2}\cap\{q_1\}^{=\pi}$, then $z\in
\{q_2,q_3\}^{=\frac\pi2}$ by Lemma \ref{lem1.2}. Now let $z$ be an
arbitrary point in
$\{q_1,q_2,q_3\}^{\geq\frac\pi2}\cap\{q_1\}^{<\pi}$. For any
geodesic $[q_1z]$, by Lemma \ref{lem1.1} we have
$|\uparrow_{q_1}^{z}\uparrow_{q_1}^{q_2}|\geq\frac\pi2$ and
$|\uparrow_{q_1}^{z}\uparrow_{q_1}^{q_3}|\geq\frac\pi2$. Together
with the condition $|\uparrow_{q_1}^{q_2}\uparrow_{q_1}^{q_3}|=\pi$,
this implies that
$$|\uparrow_{q_1}^{z}\uparrow_{q_1}^{q_2}|=|\uparrow_{q_1}^{z}\uparrow_{q_1}^{q_3}|=\frac\pi2.$$
Then by applying Theorem \ref{tct} on any triangle
$\triangle q_1q_iz$ ($i=2,3$), it is not hard to see that $|q_iz|=\frac\pi2$
(i.e. $z\in\{q_2,q_3\}^{=\frac\pi2}$);
and then by Theorem \ref{tct=}, there is a geodesic $[q_iz]$ such
that the triangle $\triangle q_1q_iz$ formed by $[q_1q_i], [q_1z]$
and $[q_iz]$ is isometric to its comparison triangle, which implies
that $\frac\pi2\leq|\uparrow_{q_i}^{z}\uparrow_{q_i}^{q_1}|<\pi$.
\hfill$\Box$

\vskip2mm

\noindent{\bf Proving $k\leq 4$ when $n=2$ in (iv) of Theorem \ref{thm2.1}}:

We consider $\Sigma_{q_1}M\in{\cal A}^1(1)$. By Lemma \ref{lem1.1},
any $\{\uparrow_{q_1}^{q_2},\cdots,\uparrow_{q_1}^{q_{k}}\}$ is a
$\frac\pi2$-separated subset in $\Sigma_{q_1}M$. If $k>4$, then by
Theorem \ref{thm0.2} we have $k=5$ and each $\uparrow_{q_1}^{q_i}$ has
an opposite direction $\uparrow_{q_1}^{q_j}$ (i.e.
$|\uparrow_{q_1}^{q_i}\uparrow_{q_1}^{q_j}|=\pi$), where $2\le i,j
\le 5$. Due to this, we can select geodesics $[q_1q_2]$ and
$[q_1q_3]$ such that
$|\uparrow_{q_1}^{q_2}\uparrow_{q_1}^{q_3}|=\pi$; so by Lemma
\ref{lem2.6}, we can select geodesics
$[q_2q_i]$ for $i=4$ and 5 such that $\frac\pi2\le
|\uparrow_{q_2}^{q_i}\uparrow_{q_2}^{q_1}|<\pi$.
Then we select an arbitrary geodesic $[q_2q_3]$, and consider
$B\triangleq\{\uparrow_{q_2}^{q_1},\uparrow_{q_2}^{q_3},\uparrow_{q_2}^{q_4},\uparrow_{q_2}^{q_{5}}\}$.
Similarly, each element of $B$ has an opposite direction in $B$. It
then follows that $|\uparrow_{q_2}^{q_1}\uparrow_{q_2}^{q_3}|=\pi$,
and similarly we can conclude that
$|\uparrow_{q_3}^{q_1}\uparrow_{q_3}^{q_2}|=\pi$. Now we let ${\cal
S}\triangleq [q_1q_2]\cup[q_2q_3]\cup[q_3q_1]$. By Lemma
\ref{lem2.2}, ${\cal S}^{\geq \frac{\pi}{2}}$ is convex in $M$; and
by Lemma \ref{lem2.4}, ${\cal S}^{\geq \frac{\pi}{2}}={\cal
S}^{=\frac{\pi}{2}}$. Then by Lemma \ref{lem2.5}, $\dim({\cal
S}^{=\frac{\pi}{2}})=0$. Note that $q_4,q_5\in {\cal S}^{\geq
\frac{\pi}{2}}$ because $q_4,q_5\in \{q_1,q_2,q_3\}^{\geq
\frac{\pi}{2}}$. Since ${\cal S}^{\geq \frac{\pi}{2}}$ is convex in
$M$, it has to hold that ${\cal S}^{\geq \frac{\pi}{2}}=\{q_4,q_5\}$
with $|q_4q_5|=\pi$ which contradicts $|q_4q_5|<\pi$. I.e., we can
conclude that $k\le4$. \hfill$\Box$

\vskip2mm

\noindent{\bf Proving $k\leq 6$ when $n=3$ in (iv) of Theorem \ref{thm2.1}}:

Similarly, if $k>6$, by Lemma \ref{lem1.1} any
$\{\uparrow_{q_1}^{q_2},$ $\cdots,\uparrow_{q_1}^{q_{7}}\}$ is a
$\frac\pi2$-separated subset in $\Sigma_{q_1}M$; and each
$\uparrow_{q_1}^{q_i}$ has an opposite direction
$\uparrow_{q_1}^{q_j}$, where $2\le i,j \le 7$. And we can find
${\cal S}\triangleq [q_1q_2]\cup[q_2q_3]\cup[q_3q_1]$ with
$|\uparrow_{q_1}^{q_2}\uparrow_{q_1}^{q_3}|=\pi$,
$|\uparrow_{q_2}^{q_1}\uparrow_{q_2}^{q_3}|=\pi$ and
$|\uparrow_{q_3}^{q_1}\uparrow_{q_3}^{q_2}|=\pi$. By Lemma
\ref{lem2.2},  ${\cal S}^{\geq \frac{\pi}{2}}$ is convex in $M$,
which implies that ${\cal S}^{\geq \frac{\pi}{2}}\in {\cal A}^m(1)$
for some $m$; and by Lemma \ref{lem2.4}, ${\cal S}^{\geq
\frac{\pi}{2}}={\cal S}^{=\frac{\pi}{2}}$. Then by Lemma
\ref{lem2.5}, we have $m\leq 1$. Note that $q_4,q_5,q_6, q_7$
belongs to ${\cal S}^{\geq \frac{\pi}{2}}$ with
$\frac\pi2\leq|q_iq_j|<\pi$ for any $4\leq i\neq j\leq 7$, which contradicts
$\dim({\cal S}^{\geq\frac{\pi}{2}})\leq 1$.
\hfill$\Box$

\subsection{\bf Estimating $k$ for $n\geq4$ in (iv) of Theorem
\ref{thm2.1}}

When we estimate $k$ for larger $n$ in (iv) of Theorem \ref{thm2.1},
the arguments in the proofs for $n=2, 3$ fail. This is because some
$\uparrow_{q_1}^{q_i}\in
\{\uparrow_{q_1}^{q_2},\cdots,\uparrow_{q_1}^{q_{k}}\}$ maybe have
no opposite direction in
$\{\uparrow_{q_1}^{q_2},\cdots,\uparrow_{q_1}^{q_{k}}\}$ for larger
$n$. In order to overcome this, we need the following proposition.

\begin{lemma}\label{lem2.7}
Let $N$ be a convex subset in $M\in{\cal A}^n(1)$ with
$\dim(N)=n-1$. If $N^{=\frac\pi2}\neq\emptyset$, then
$N^{=\frac\pi2}=\{p\}$ or $\{p_1,p_2\}$, and

\noindent{\rm (i)}  if $N^{=\frac\pi2}=\{p\}$ (resp. $\{p_1,p_2\}$),
then there are at most two (resp. a unique) geodesics between $p$
(resp. $p_i$) and any interior point $x$ of $N$;

\noindent{\rm(ii)}  if $N^{=\frac\pi2}=\{p_1,p_2\}$, and if $N$ is
complete and $N$ has empty boundary, then $M=\{p_1,p_2\}*N$;

\noindent{\rm(iii)}  if $N^{=\frac\pi2}=\{p_1,p_2\}$, then
$\{p_1,p_2\}^{\geq\frac\pi2}=\{p_1,p_2\}^{=\frac\pi2}$.
\end{lemma}

\noindent{\it Proof}. We will prove this by the induction on $n$.
Obviously, the lemma is true if $n=1$, so we assume $n>1$.

Note that by the first variation formula ([BGP]), for any $y\in
N^{=\frac\pi2}$, $x\in N^\circ$ and any geodesic $[yx]$, we have
$$|\uparrow_{x}^y\xi|\ge\frac\pi2,\ \forall\ \xi\in\Sigma_xN.$$
Note that $\Sigma_xN$ is convex in $\Sigma_xM$ because $N$ is convex
in $M$, and $\Sigma_xN$ has empty boundary because $x\in N^\circ$
([BGP]). Then by Lemma \ref{lem2.3}, $|\uparrow_{x}^y\xi|=\frac\pi2$
in fact, i.e., $\uparrow_{x}^y\in (\Sigma_xN)^{=\frac\pi2}$; and by
the inductive assumption we know that $(\Sigma_xN)^{=\frac\pi2}$
contains at most two points in $\Sigma_xM$. Hence, we can conclude
that $N^{=\frac\pi2}=\{p\}$ or $\{p_1,p_2\}$, and that (i) holds.

(ii) Due to (i), this is a special case of Proposition \ref{prop2.8} below.

(iii) Given any fixed $x\in N^\circ$, from the above, we know that
$(\Sigma_xN)^{=\frac\pi2}=\{\uparrow_{x}^{p_1},\uparrow_{x}^{p_2}\}$.
Note that $\Sigma_xN$ is convex and complete in $\Sigma_x M$ and has
empty boundary. Then by (ii), we have
$\Sigma_xM=\{\uparrow_{x}^{p_1}, \uparrow_{x}^{p_2}\}*\Sigma_x N$.
Hence, for any $z\in \{p_1,p_2\}^{\geq\frac\pi2}$ and a geodesic
$[xz]$, without loss of generality we can assume that
$|\uparrow_{x}^{p_1}\uparrow_{x}^{z}|\leq \frac\pi2$. On the other
hand, by applying Theorem \ref{tct} on $\triangle p_1xz$, we have
$|\uparrow_{x}^{p_1}\uparrow_{x}^{z}|\geq \frac\pi2$ (note that
$|p_1z|\geq\frac\pi2$ and $|p_1x|=\frac\pi2$). It then follows that
$|\uparrow_{x}^{p_1}\uparrow_{x}^{z}|=\frac\pi2$, which implies that
$|\uparrow_{x}^{p_2}\uparrow_{x}^{z}|=\frac\pi2$ too. Again by
applying Theorem \ref{tct} on $\triangle p_ixz$ ($i=1,2$), we
conclude that $|p_iz|=\frac\pi2$, i.e.,
$\{p_1,p_2\}^{\geq\frac\pi2}=\{p_1,p_2\}^{=\frac\pi2}$. \hfill$\Box$

\begin{prop}\label{prop2.8} Let $X, Y$ be two complete convex subsets in $M\in{\cal A}^n(1)$.
Then $M=X*Y$ if\ \ {\rm (i)} $\dim(X)+\dim(Y)+1=n$;

\noindent{\rm (ii)} $X$ and $Y$ have empty boundary;

\noindent{\rm (iii)} $|xy|=\frac\pi2$ for any $x\in X$ and $y\in Y$;

\noindent{\rm (iv)} there is a unique geodesic between any $x\in X$
and $y\in Y$.
\end{prop}

Proposition \ref{prop2.8} is due to the definition of the metric of
the join ([BGP]). For a detailed proof, one can refer to [SSW].

\vskip1mm

Now we prove that, in (iv) of Theorem
\ref{thm2.1}, $k$ has the desired upper bound for $n\geq4$.

\vskip2mm

\noindent{\bf Proof for the estimate of $k$ for $n\geq4$ in (iv) of
Theorem \ref{thm2.1}}:

We consider $\Sigma_{q_1}M\in{\cal A}^{n-1}(1)$. By Lemma
\ref{lem1.1}, any
$\{\uparrow_{q_1}^{q_2},\cdots,\uparrow_{q_1}^{q_k}\}$
is a $\frac\pi2$-separated subset in $\Sigma_{q_1}M$. If
$\{\uparrow_{q_1}^{q_2},\cdots,\uparrow_{q_1}^{q_{k}}\}$
has no two opposite directions
(i.e. $|\uparrow_{q_1}^{q_i}\uparrow_{q_1}^{q_{j}}|<\pi$ for any $2\leq i\neq j\leq k$),
then by induction we have that
$k-1\leq 3(l-1)+1$ (resp. $3l$) for $n-1=2l-2$ (resp. $2l-1$), i.e.,
$k\leq 3l-1$ (resp. $3l+1$) for $n=2l-1$ (resp. $2l$). Hence,
$k\leq 3l$ (resp. $3l+1$) for $n=2l-1$ (resp. $2l$).

From the above, we conclude that if $k>3l$ (resp. $3l+1$) for
$n=2l-1$ (resp. $2l$) (Hint: In fact, ``$k>3l$'' can be changed to
``$k\geq3l$'' when $n=2l-1$), then
\begin{equation}\label{eqn2.9}
\text{\it for any $1\leq i\leq k$, any $\{\uparrow_{q_i}^{q_j}|
1\leq  j\leq k, j\ne i \}$ has two opposite directions.}
\end{equation}

Now we assume that $k>3l$ (resp. $3l+1$) for $n=2l-1$ (resp. $2l$),
and we consider
$\{\uparrow_{q_{1}}^{q_2},\cdots,\uparrow_{q_1}^{q_{k}}\}$. By
(\ref{eqn2.9}), without loss of generality, we can assume that
$$|\uparrow_{q_{1}}^{q_2}\uparrow_{q_{1}}^{q_3}|=\pi.$$ By
Lemma \ref{lem2.2}, $X\triangleq\{q_1,q_2,q_3\}^{\geq\frac\pi2}$ is
convex in $M$, so $X\in {\cal A}^m(1)$ for some $m$. And by Lemma
\ref{lem2.6}, $X\subseteq \{q_{2},q_{3}\}^{=\frac\pi2}$, so $m\leq
n-1$ by Lemma \ref{lem2.5}. Note that $\{q_4,\cdots,q_k\}$ is a
$\frac\pi2$-separated subset in $X$.

If $m\leq n-2$, then $m\leq 2l-3$ (resp. $2l-2$) for
$n=2l-1$ (resp. $2l$). By the inductive assumption on $X$, we
conclude that $k-3\leq 3(l-1)$ (resp. $3(l-1)+1$) for $m\leq
2l-3$ (resp. $2l-2$), which contradicts the assumption ``$k>3l$ (resp. $3l+1$)''.

If $m=n-1$,  then by Lemma \ref{lem2.7} it has to hold that
$X^{=\frac\pi2}=\{q_2,q_3\}$. And by (iii) and (i) of Lemma
\ref{lem2.7}, we have
$Y\triangleq\{q_{2},q_{3}\}^{\ge\frac\pi2}=\{q_{2},q_{3}\}^{=\frac\pi2}$,
and there is a unique geodesic between $q_{2}$ and any interior
point of $Y$. This implies that $\{q_{2}\}*Y$ can be isometrically
embedded into $M$ (note that $Y$ is convex in $M$ by Lemma
\ref{lem2.2}, and see Proposition \ref{prop2.8} or refer to [SSW]).
Note that $\{q_1,q_4,\cdots,q_k\}\subset Y$. Then we can find a
$\frac\pi2$-separated subset
$\{\uparrow_{q_{2}}^{q_1},\uparrow_{q_{2}}^{q_{4}},\cdots,\uparrow_{q_{2}}^{q_{k}}\}$
in $\Sigma_{q_{2}}M$ with
$|\uparrow_{q_{2}}^{q_{i}}\uparrow_{q_{2}}^{q_{j}}|=|q_iq_j|$ for
any $i,j\in\{1,4,\cdots,k\}$. Since $|q_iq_j|<\pi$, due to
(\ref{eqn2.9}) there is a geodesic $[q_{2}q_{3}]$ such that
$|\uparrow_{q_{2}}^{q_{i_0}}\uparrow_{q_{2}}^{q_{3}}|=\pi$ for some
$i_0\in \{1,4,\cdots,k\}$. By Lemma \ref{lem2.6}, $Z\triangleq
Y\cap\{q_{i_0}\}^{\geq\frac\pi2}$ belongs to
$Y\cap\{q_{i_0}\}^{=\frac\pi2}$, which implies that $\dim(Z)\leq
m-1=n-2$ by Lemma \ref{lem2.5} (note that $Z$ is convex in $Y$ (and
$M$)). By the inductive assumption on $Z$ (note that
$\{q_1,\cdots,q_k\}\setminus\{q_2,q_3,q_{i_0}\}$ is a
$\frac\pi2$-separated subset of $Z$), we conclude that $k-3\leq
3(l-1)$ (resp. $3(l-1)+1$) for $n-2=2l-3$ (resp. $2l-2$), which
contradicts the assumption ``$k>3l$ (resp. $3l+1$)''.

Since all contradictions are gotten under the assumption ``$k>3l$
(resp. $3l+1$)'' for $n=2l-1$ (resp. $2l$), we conclude that
$k\leq3l$ (resp. $3l+1$). \hfill $\Box$

\begin{remark}\label{rem2.10}{\rm In proving that $k$ has
the  desired upper bound in (iii) of Theorem \ref{thm2.1}, a main
difference to the above proof is that we will consider
$\{\uparrow_{q_{h+1}}^{q_1},\cdots,\uparrow_{q_{h+1}}^{q_h},
\uparrow_{q_{h+1}}^{q_{h+2}},\cdots,$ $\uparrow_{q_{h+1}}^{q_{k}}\}$
in $\Sigma_{q_{h+1}}M$ (here it will be better to point out that,
by Lemma \ref{lem1.2}, it is not hard to see that we have
$i, j>h$ if $|\uparrow_{q_{h+1}}^{q_{i}}\uparrow_{q_{h+1}}^{q_{j}}|=\pi$).}
\end{remark}

\subsection{Proof for the structure classification in (iv) of Theorem \ref{thm2.1}}

\noindent{\bf Proof for $n=2$ and $k=4$}:

We only need to prove that there is a point in $\{q_1,\cdots,q_4\}$,
say $q_4$, such that $M=\{q_4\}*S^1$ with $\{q_1,q_2,q_3\}\subset
S^1$ if $M$ has nonempty boundary.

{\bf Claim 1}: {\it There is a $q_i$ with $q_i\in\partial M$}. If
$|q_1\partial M|\geq\frac\pi2$, then it is easy to see that, by
Lemma 1.4, $q_i\in\partial M$ for $2\leq i\leq 4$. If $|q_1\partial
M|<\frac\pi2$, we select $p\in \partial M$ such that $|q_1p|=|q_1\partial M|$.
Like Proof \ref{pf1.3}, we can prove that
$\Sigma_pM=\{\uparrow_{p}^{q_1}\}*\partial(\Sigma_pM)$, and that
$|pq_j|\geq\frac{\pi}{2}$ for $j\neq1$. If $|pq_j|=\pi$ for some
$j\neq1$, then $M=\{p,q_j\}*A$ for some arc $A$ (Lemma \ref{lem1.2}) which implies
$q_j\in\partial M$. If $|pq_j|<\pi$ for $j=2,3$ and 4, then
by Lemma \ref{lem1.1}, any $\{\uparrow_{p}^{q_2},
\uparrow_{p}^{q_3},\uparrow_{p}^{q_4}\}$ is a $\frac\pi2$-separated
subset in $\Sigma_pM (=\{\uparrow_{p}^{q_1}\}*\partial(\Sigma_pM))\in {\cal A}^1(1)$. It therefore
follows that, without loss of generality, we can assume that
$\uparrow_{p}^{q_2}$ and $\uparrow_{p}^{q_3}$ belong to
$\partial(\Sigma_pM)$ with
$|\uparrow_{p}^{q_2}\uparrow_{p}^{q_3}|=\pi$. Note that
$\uparrow_{p}^{q_2}\in\partial(\Sigma_pM)$ implies that there exists
a geodesic $[pq_2]$ belonging to $\partial M$ ([BGP]), i.e., Claim 1
is verified.

{\bf Claim 2}: {\it There is a $q_i$ with $q_i\not\in\partial M$}.
If this is not true, we select any $q_i$, say $q_2$, and consider $\Sigma_{q_2}M\in {\cal
A}^1(1)$ which has nonempty boundary. By Lemma
\ref{lem1.1}, any $\{\uparrow_{q_2}^{q_1},
\uparrow_{q_2}^{q_3},\uparrow_{q_2}^{q_4}\}$ is a
$\frac\pi2$-separated subset in $\Sigma_{q_2}M$; so without loss of
generality we can assume that $\uparrow_{q_2}^{q_1}$ and
$\uparrow_{q_2}^{q_3}$ belong to $\partial(\Sigma_{q_2}M)$ with
$|\uparrow_{q_2}^{q_1}\uparrow_{q_2}^{q_3}|=\pi$. This implies that
there exists a unique geodesic between $q_2$ and $q_k$ ($k=1,3$) and $[q_2q_k]$
belongs to $\partial M$ (Hint: if $|pq|=\pi$ for $p,q\in M\in {\cal
A}^n(1)$, then $q$ is the unique point in $M$ such that $|pq|=\pi$
(Lemma \ref{lem1.2})). Hence, if Claim 2 is not true, then we can rearrange
$q_1,\cdots,q_4$ such that there is a unique geodesic
$[q_{i}q_{i+1}]$ between $q_{i}$ and $q_{i+1}$ for $i=1,\cdots,4$
(where $q_5=q_1$) with $\cup_{i=1}^4[q_{i}q_{i+1}]\subseteq\partial
M$ and $|\uparrow_{q_i}^{q_{i-1}}\uparrow_{q_i}^{q_{i+1}}|=\pi$
(where $q_{0}=q_4$). By Lemma \ref{lem2.6} and its proof, we have
$|\uparrow_{q_1}^{q_3}\uparrow_{q_1}^{q_2}|=
|\uparrow_{q_1}^{q_3}\uparrow_{q_1}^{q_4}|=
|\uparrow_{q_3}^{q_1}\uparrow_{q_3}^{q_2}|=\frac\pi2$ and
$|q_3q_2|=\frac\pi2$. Then by Theorem \ref{tct=}, there is a geodesic $[q_1q_3]$
such that the triangle $\triangle q_2q_1q_3$ (formed by $[q_1q_2]$, $[q_2q_3]$ and $[q_1q_3]$)
is isometric to its comparison triangle,
so $|\uparrow_{q_2}^{q_1}\uparrow_{q_2}^{q_3}|=|q_1q_3|<\pi$ which contradicts
$|\uparrow_{q_2}^{q_1}\uparrow_{q_2}^{q_3}|=\pi$. Hence, Claim 2 has
to hold.

Due to Claims 1 and 2, we can assume that $q_4\not\in\partial M$ and $q_2\in \partial M$.
Then by the proof of Claim 2, we can conclude that $q_1,q_3\in\partial M$ too, and there are geodesics $[q_{1}q_{2}],
[q_2q_3]$ and $[q_3q_1]$ such that ${\cal
S}\triangleq[q_{1}q_{2}]\cup[q_2q_3]\cup[q_3q_1]\subseteq\partial M$
and $|\uparrow_{q_1}^{q_{3}}\uparrow_{q_1}^{q_2}|=
|\uparrow_{q_2}^{q_1}\uparrow_{q_2}^{q_3}|=
|\uparrow_{q_3}^{q_{2}}\uparrow_{q_3}^{q_1}|=\pi$. By Lemma
\ref{lem2.4}, we have ${\cal S}^{\geq\frac\pi2}={\cal
S}^{=\frac\pi2}$. Since $|q_4q_i|\geq\frac\pi2$ for $i=1,2,3$, by
Theorem \ref{tct} we have $q_4\in{\cal S}^{\geq\frac\pi2}$, and so
$|q_4p|=\frac\pi2$ for all $p\in{\cal S}$. By the first variation
formula ([BGP]) together with ${\cal S}\subseteq\partial M$, we have
that any geodesic $[q_4p]$ is perpendicular to ${\cal S}$ at any
$p\in {\cal S}$, and that there is a unique direction which is
perpendicular to ${\cal S}$ at any $p\in {\cal S}$. This implies
that there is a unique geodesic between $q_4$ and any $p\in{\cal
S}$, and that any point in $M$ lies in a geodesic $[q_4p]$ for some
$p\in {\cal S}$. Then by Theorem \ref{tct=}, we can conclude that
$M=\{q_4\}*{\cal S}$, i.e., $M=\{q_4\}*S^1$ with
$\{q_1,q_2,q_3\}\subset S^1$. \hfill$\Box$

\vskip2mm

In the proof for $n\geq3$, we need a more technical result than Proposition
\ref{prop2.8}.

\begin{lemma}\label{lem2.11}
Let $X$ and a circle $S^1$ be two convex subsets in $M\in{\cal
A}^n(1)$. Then $M=X*S^1$ if\ \ {\rm (i)} $\dim(X)=n-2$; {\rm (ii)}
$X$ is complete and has empty boundary; {\rm (iii)}  the perimeter
of $S^1$ is bigger than $\pi$; {\rm (iv)} $|xy|=\frac\pi2$ for any
$x\in X$ and $y\in S^1$.
\end{lemma}

In Lemma \ref{lem2.11}, if condition (iii) is canceled, Rong-Wang
proved that there are $\hat X\in{\cal A}^{n-2}(1)$, $\hat
S^1\in{\cal A}^{1}(1)$ and a cyclic group $\Gamma$ which acts by
isometries on $\hat X$ and $\hat S^1$ such that $X=\hat X/\Gamma$,
$S^1=\hat S^1/\Gamma$ and $M=(\hat X*\hat S^1)/\Gamma$ (cf. [RW]).
Note that this implies Lemma \ref{lem2.11}. However, for the
convenience of readers, we will give the proof of Lemma
\ref{lem2.11} in Appendix.

Now we prove the latter part of (iv) for $n\geq3$ in Theorem
\ref{thm2.1}.

\vskip2mm

\noindent{\bf Proof for $n\geq3$}:

\vskip1mm

\noindent{Case 1}: $n=2l-1$  and $k=3l$ with $l\geq 2$.

Note that (\ref{eqn2.9}) still holds when $n=2l-1$ and $k=3l$ (see
the hint before (\ref{eqn2.9})). Like the proof for the estimate of
$k$ for $n\geq4$ in (iv) of Theorem \ref{thm2.1}, we can select
$q_1,q_2,q_3$ such that
$X\triangleq\{q_{1},q_{2},q_{3}\}^{\geq\frac\pi2}\in {\cal A}^m(1)$
belongs to $\{q_{2},q_{3}\}^{=\frac\pi2}$. Since $X$ is convex in
$M$ (Lemma \ref{lem2.2}), we have that $m\leq n-1$ (Lemma 2.5). Note
that $\{q_{4},\cdots,q_k\}$ is a $\frac\pi2$-separated subset in $X$
with $|q_iq_j|<\pi$ ($i\neq j$), which implies that $m\geq n-2$ (by
the estimate of $k$ in (iv) of Theorem \ref{thm2.1}). It then
follows that $m=n-2$ or $n-1$.

If $m=n-2$, then by induction we conclude that $X$ is isometric
to $S^1_2*\cdots*S^1_l$, where $S^1_j$ ($2\leq j\leq l$) has
perimeter $\geq\frac{3\pi}{2}$ (and we can rearrange
$q_4,\cdots,q_{3l}$ such that $q_{3j-2}, q_{3j-1}, q_{3j}\in
S_j^1$). Of course, $X$ has empty boundary. Since $X$ is convex in
$M$, by Lemmas \ref{lem2.3} and \ref{lem2.2},
$X^{\geq\frac\pi2}=X^{=\frac\pi2}$ and $X^{\geq\frac\pi2}$ is convex
in $M$, and so by Lemma \ref{lem2.5} $\dim(X^{=\frac\pi2})\leq 1$.
Obviously, $\{q_{1},q_{2},q_{3}\}\subset X^{=\frac\pi2}$, so it has
to hold that $X^{=\frac\pi2}$ is a circle with perimeter
$\geq\frac{3\pi}{2}$, denoted by $S^1_1$. Then by Lemma
\ref{lem2.11}, it follows that $M$ is isometric to
$S^1_1*\cdots*S^1_{l}$ with $q_{3j-2}, q_{3j-1}, q_{3j}\in S_j^1$
for $1\leq j\leq l$.

If $m=n-1$, like the case ``$m=n-1$'' in the proof for the estimate
of $k$ for $n\geq4$ in (iv) of Theorem \ref{thm2.1}, we have that
$X^{=\frac\pi2}=\{q_2,q_3\}$ and
$Y\triangleq\{q_{2},q_{3}\}^{\ge\frac\pi2}=\{q_{2},q_{3}\}^{=\frac\pi2}$.
And we can select $q_{i_0}$ ($i_0\neq 2,3$) such that $Z\triangleq
Y\cap\{q_{i_0}\}^{\geq\frac\pi2}=Y\cap\{q_{i_0}\}^{=\frac\pi2}$; and
$Z$ is convex in $M$ with $\dim(Z)\leq m-1=n-2$; and
$\{q_1,\cdots,q_{3l}\}\setminus\{q_2,q_3,q_{i_0}\}$ is a
$\frac\pi2$-separated subset of $Z$.  By induction, $Z$ is isometric to $S^1_2*\cdots*S^1_l$, and
similarly we can derive that $M$ is isometric to
$S^1_1*\cdots*S^1_{l}$ with $S^1_j$ having perimeter
$\geq\frac{3\pi}{2}$ (and we can rearrange $q_1,\cdots,q_{3l}$ such
that $q_{3j-2}, q_{3j-1}, q_{3j}\in S_j^1$ for $1\leq j\leq l$).

\vskip1mm

\noindent{Case 2}: $n=2l$ and $k=3l+1$ with $l\geq 2$.

\vskip1mm

In this case, (\ref{eqn2.9}) may not hold. We will give discussions according
to ``(\ref{eqn2.9}) holds'' and ``(\ref{eqn2.9}) does not hold''.

Subcase 1: (\ref{eqn2.9}) holds.

Like in Case 1 ($n=2l-1$ and $k=3l$), we can find convex and
complete subset $X$ or $Z\in{\cal A}^{n-2}(1)$ in which
$\{q_4,\cdots,q_{3l+1}\}$ or $\{q_1,\cdots,q_{3l+1}\}\setminus
\{q_{i_0},q_{2},q_{3}\}$ are $\frac\pi2$-separated subsets
respectively. Without loss of generality, we assume that such an $X$
is found.

If $X$ has empty boundary, then $X^{\geq\frac\pi2}=X^{=\frac\pi2}$
(Lemma \ref{lem2.3}) which is convex in $M$ (Lemma \ref{lem2.2}),
and thus $\dim(X^{\geq\frac\pi2})\leq 1$ (Lemma \ref{lem2.5}). Note
that $\{q_{1},q_{2},q_{3}\}$ is a $\frac\pi2$-separated subset of
$X^{\geq\frac\pi2}$, so it has to hold that $X^{\geq\frac\pi2}$ is a
circle $S^1$ with perimeter $\geq\frac{3\pi}{2}$. Hence, by Lemma
\ref{lem2.11} we have $M=X*S^1$. This implies that $M$ has the
desired structure because $X\in {\cal A}^{n-2}(1)$ has the desired
structure by induction.

If $X$ has nonempty boundary, then by induction we can rearrange
$q_4,\cdots,q_{3l+1}$ such that $X$ is isometric to
$\{q_{3l+1}\}*S^1_2*\cdots*S^1_{l}$ with $S^1_j$ having perimeter
$\geq\frac{3\pi}{2}$ and $q_{3j-2}, q_{3j-1}, q_{3j}\in S_j^1$ for
$2\leq j\leq l$. Now we consider $(S^1_2)^{\geq\frac\pi2}$. Note
that $S^1_2$ is convex in $M$ (because $X$
($=\{q_{3l+1}\}*S^1_2*\cdots*S^1_{l}$) is convex in $M$). Then
$W\triangleq(S^1_2)^{\geq\frac\pi2}$ is convex in $M$ (Lemma
\ref{lem2.2}), and $W=(S^1_2)^{=\frac\pi2}$ (Lemma \ref{lem2.3}),
and $\dim(W)\leq n-2$ (Lemma \ref{lem2.5}). Note that
$\{q_1,q_2,q_3,q_7,\cdots,q_{3l+1}\}$ is a $\frac\pi2$-separated
subset of $W$, so $\dim(W)=n-2$ and $W\in {\cal A}^{n-2}(1)$ has the
desired structure by induction. If $W$ has empty boundary, then by
Lemma \ref{lem2.11} $M$ is isometric to $S^1_2*W$, and thus $M$ has
the desired structure. If $W$ has nonempty boundary, then by
induction we can assume that
$W=\{q_{3l+1}\}*S^1_1*S^1_3*\cdots*S^1_{l}$ with $q_1, q_2, q_3\in
S_1^1$, and from the proof of Lemma \ref{lem2.11} (see Remark A.1 in
Appendix) $S^1_2*W$ can be isometrically embedded into $M$. And it
is not hard to see that $M=\{q_{3l+1}\}*S^1_1*S^1_2*\cdots*S^1_{l}$
if $M$ has nonempty boundary.

\vskip1mm

Subcase 2: (\ref{eqn2.9}) does not hold.

In the proof, we will need the following lemma.

\begin{lemma}\label{lem2.12}
Let $X*Y$ be a join with $X, Y\in \mathcal{A}(1)$. If
$Q\triangleq\{q_1,\cdots,q_k\}$ is a $\frac\pi2$-net of $X$ (i.e.
for any $x\in X$ there is $q_i\in Q$ such that $|xq_i|<\frac\pi2$),
then $Q^{\geq\frac\pi2}= Y$.
\end{lemma}

\noindent{Proof}. It sufficed to show that $Q^{\geq\frac\pi2}\subset Y$. Let $p$ be any  point in $Q^{\geq\frac\pi2}$. If
$p\not\in Y$, then $p=[(x,y,t)]$ with $x\in X$, $y\in Y$ and
$t<\frac\pi2$. Select $q_i\in Q$ such that $|xq_i|<\frac\pi2$ and a
geodesic $[xq_i]$. By the definition of the metric of the join, we
know that $|\uparrow_x^{q_i}\uparrow_x^{p}|=\frac\pi2$ (cf. [SSW]).
Then by Theorem \ref{tct} on the triangle $\triangle pxq_i$, we conclude that
$|pq_i|<\frac\pi2$ (note that $|xp|=t<\frac\pi2$ and
$|xq_i|<\frac\pi2$) which contradicts $p\in Q^{\geq\frac\pi2}$.
Hence, it has to hold that $p\in Y$. \hskip1cm$\blacksquare$

\vskip1mm

Now we continue the proof under the assumption ``(\ref{eqn2.9}) does
not hold''.

Since (\ref{eqn2.9}) does not hold, without loss of generality, we
can select
$\{\uparrow_{q_{3l+1}}^{q_1},\cdots,\uparrow_{q_{3l+1}}^{q_{3l}}\}$
in $\Sigma_{q_{3l+1}}M\in{\cal A}^{n-1}(1)$ such that
$|\uparrow_{q_{3l+1}}^{q_i}\uparrow_{q_{3l+1}}^{q_j}|<\pi$ for any
$1\leq i\neq j\leq 3l$. On the other hand,
$\{\uparrow_{q_{3l+1}}^{q_1},\cdots,\uparrow_{q_{3l+1}}^{q_{3l}}\}$
is a $\frac{\pi}{2}$-separated subset of $\Sigma_{q_{3l+1}}M\in{\cal
A}^{n-1}(1)$ (Lemma \ref{lem1.1}). Then by the case ``$n=2l-1$ and
$k=3l$'' (here $n-1=2l-1$), we can conclude that
$\Sigma_{q_{3l+1}}M=\bar S^1_1*\cdots*\bar S^1_{l}$ with $\bar
S^1_i$ having perimeter $\geq\frac{3\pi}{2}$ and
$\uparrow_{q_{3l+1}}^{q_{3i-2}},\uparrow_{q_{3l+1}}^{q_{3i-1}},\uparrow_{q_{3l+1}}^{q_{3i}}\in
\bar S_i^1$.

Let $A\triangleq\{q_{3l+1},q_{3l}, q_{3l-1}, q_{3l-2}\}^{\ge
\frac\pi2}$, which is convex in $M$ (Lemma \ref{lem2.2}). If there
is $p\in A$ such that $|q_{3l+1}p|=\pi$, then $M=\{q_{3l+1},p\}*L$
for some $L\in{\cal A}^{n-1}(1)$ (Lemma \ref{lem1.2}). Note that we can assume that
$\{q_1,\cdots,q_{3l}\}$ is a $\frac{\pi}{2}$-separated subset of $L$
(otherwise we can replace $q_i$ with the the point $[pq_i]\cap L$ for $i=1,2,\cdots,3l$).
Similarly, by the case ``$n=2l-1$ and $k=3l$'', we can conclude that
$L=S^1_1*\cdots*S^1_{l}$ with $S^1_i$ having perimeter
$\geq\frac{3\pi}{2}$ and $q_{3i-2},q_{3i-1},q_{3i}\in S_i^1$, and so
$\{q_{3l+1}\}*S^1_1*\cdots*S^1_{l}$ can be isometrically embedded
into $M$. Now we assume that $|q_{3l+1}p|<\pi$ for any $p\in A$. By
Lemma \ref{lem1.1}, for any geodesic $[q_{3l+1}p]$, we have
$|\uparrow_{q_{3l+1}}^{p}\uparrow_{q_{3l+1}}^{q_j}|\geq\frac\pi2$
for $j=3l-2, 3l-1, 3l$. Note that
$\{\uparrow_{q_{3l+1}}^{q_{3l-2}},\uparrow_{q_{3l+1}}^{q_{3l-1}},\uparrow_{q_{3l+1}}^{q_{3l}}\}$
is a $\frac\pi2$-net of $\bar S_l^1$, so by Lemma \ref{lem2.12}
$$\uparrow_{q_{3l+1}}^{p}\in \bar S^1_1*\cdots*\bar S^1_{l-1}.$$
Next we consider $B\triangleq(A\cup\{q_{3l+1}\})^{\ge \frac\pi2}$,
which is also convex in $M$. Similarly, we assume that
$|q_{3l+1}r|<\pi$ for any $r\in B$; so for any geodesic $[q_{3l+1}r]$,
$|\uparrow_{q_{3l+1}}^{r}\uparrow_{q_{3l+1}}^{q_{j}}|\geq\frac\pi2$
for any $1\leq j\leq 3l-3$. And it is not hard to see that
$\{\uparrow_{q_{3l+1}}^{q_{1}},\cdots,\uparrow_{q_{3l+1}}^{q_{3l-3}}\}$
is a $\frac\pi2$-net of $\bar S^1_1*\cdots*\bar S^1_{l-1}$, so by
Lemma \ref{lem2.12}
$$\uparrow_{q_{3l+1}}^{r}\in \bar S^1_{l}.$$
It therefore follows that
$|\uparrow_{q_{3l+1}}^{p}\uparrow_{q_{3l+1}}^{r}|=\frac\pi2$. Then
by Theorem \ref{tct} on any $\triangle q_{3l+1}pr$, we have $|pr|=\frac\pi2$ (note
that $|q_{3l+1}p|, |q_{3l+1}r|, |pr|\geq\frac\pi2$). Hence, by Lemma
\ref{lem2.5}, we have $$\dim(A)+\dim(B)\leq n-1.$$ Note that
$\{q_1,\cdots,q_{3l-3}\}\subset A$ and
$\{q_{3l-2},q_{3l-1},q_{3l}\}\subset B$, and $A$ and $B$ are convex
in $M$. This implies that $\dim(A)\geq n-3$ (see Subsection 2.3) and
$\dim(B)\geq 1$. It then has to hold that either $\dim(A)=n-3$, or
$\dim(A)=n-2$ and $\dim(B)=1$.

If $\dim(A)=n-3$, then by the case ``$n=2l-1$ and $k=3l$'' (here
$n-3=2l-3$) we get that $A=S^1_1*\cdots*S^1_{l-1}$ with $S^1_i$
having perimeter $\geq\frac{3\pi}{2}$ and
$q_{3i-2},q_{3i-1},q_{3i}\in S_i^1$. If $\dim(A)=n-2$ and
$\dim(B)=1$, then $B$ ($\ni q_{3l-2},q_{3l-1},q_{3l}$) is a circle
with perimeter $\geq\frac{3\pi}{2}$. Now we let ${\cal S}$ denote
$S_1^1$ or the circle $B$ (which is convex in $M$), and we consider
$C\triangleq ({\cal S})^{\ge \frac\pi2}$, which is convex in $M$
(Lemma \ref{lem2.2}). By Lemmas \ref{lem2.3} and \ref{lem2.5},
$C=({\cal S})^{=\frac\pi2}$ and $\dim(C)\leq n-2$. Note that
$\{q_4,\cdots,q_{3l+1}\}\subset C$ or
$\{q_1,\cdots,q_{3l-3},q_{3l+1}\}\subset C$, so we have $\dim(C)\geq
n-2$ (see Subsection 2.3). Hence, we have $\dim(C)=n-2$, and so by
Lemma \ref{lem2.11} or its proof (see Remark A.1) we get that
$M={\cal S}*C$ (if $C$ has empty boundary) or ${\cal S}*C$ can be
isometrically embedded into $M$ (if $C$ has nonempty boundary). On
the other hand, $C$ has the desired structure by induction. Hence,
it is not hard to see that $M$ has the desired structure.
\hfill$\Box$

\vskip2mm

So far we have finished the proof of Theorem \ref{thm2.1} (which
implies Theorem B). In the rest of this section, we will give the
proof of Corollary C.

\vskip1mm

\noindent{\bf Proof of Corollary C}:

Since $M\in{\cal A}^n(1)$, by Theorem B, we have $k\leq 3l\
(\text{resp. } 3l+1) \text{ for } n=2l-1\ (\text{resp. } 2l)$;
moreover, since $M$ has empty boundary (because $M$ is closed) and
$n>2$, $M$ is isometric to $S_1^1*\cdots
*S_l^1$ (resp. $S_1^1*\cdots
*S_{l-1}^1*N$ for some $N\in {\cal A}^2(1)$ or $\{q_{3l+1}\}*S_1^1*\cdots
*S_{l}^1$ can be isometrically embedded into $M$) with $l\geq2$ if $k=3l\
(\text{resp. } 3l+1) \text{ for } n=2l-1\ (\text{resp. } 2l)$. Given
a join $X*Y$ with $X, Y\in{\cal A}(1)$, from the definition of the
metric of the join (cf. [SSW]), we know that
$\Sigma_x(X*Y)=(\Sigma_xX)*Y$ for any $x\in X\subset X*Y$. Since $M$
is a Riemannian manifold, $\Sigma_pM$ is isometric to the unit
sphere $\Bbb S^{n-1}$ for any $p\in M$. Therefore, if $k=3l\
(\text{resp. } 3l+1) \text{ for } n=2l-1\ (\text{resp. } 2l)$, then
each $S^1_i$ is a great circle (i.e. having perimeter equal to
$2\pi$) which can be isometrically embedded into $M$. It then
follows from the Maximum Diameter Theorem ([CE]) that $M$ is
isometric to the unit sphere $\Bbb S^n$. \hfill$\Box$

%%%%%%%%%%%%%%%%%%%%%%%%%%%%%%%%%%%%%%% Section 3  %%%%%%%%%%%%%%%%%%%%%%%%%%%%%%%%%%%%%%%

\section{Proofs of Theorem D and Corollary E}

We first give an interesting and key lemma, which may be known to
experts.

\vskip2mm

\begin{lemma}\label{lem3.1}  Let $M\in \mathcal{A}^n(1)$, and let
$\mathbb{S}^k$ be the $k$-dimensional unit sphere. If there exists a
noncontractive map $f:\mathbb{S}^k \to M$, then $f$ is an isometrical
embedding.
\end{lemma}
\vskip2mm

\noindent{\it Proof}. It suffices to show that $|f(x)f(y)|=|xy|$ for
any $x,y\in \mathbb{S}^k$. Note that $\diam(M)\leq\pi$ because $M\in
\mathcal{A}^n(1)$ ([BGP]); so if $|xy|=\pi$, then
$|f(x)f(y)|=\pi=|xy|$ because $f$ is a noncontractive map (note that
$\mathbb{S}^0$ consists of two points with distance equal to $\pi$).
Now we assume that $|xy|<\pi$. Let $z$ be the antipodal point of the
middle point of $[xy]$ (in $\mathbb{S}^k$, $[xy]$ is the unique
geodesic between $x$ and $y$ if $|xy|<\pi$). Note that
\begin{equation}\label{eqn3.2}
|xy|+|xz|+|zy|=2\pi.
\end{equation}
Since $f$ is a noncontractive map, we have
\begin{equation}\label{eqn3.3}
|f(x)f(y)|\geq |xy|,\quad |f(x)f(z)|\geq |xz|,\quad |f(z)f(y)|\geq
|zy|,
\end{equation}
and thus
$$|f(x)f(y)|+ |f(x)f(z)|+|f(z)f(y)|\geq2\pi.$$
On the other hand, because $M\in \mathcal{A}^n(1)$, we have ([BGP])
$$|f(x)f(y)|+ |f(x)f(z)|+|f(z)f(y)|\leq 2\pi.$$
It then follows that
\begin{equation}\label{eqn3.4}
|f(x)f(y)|+ |f(x)f(z)|+|f(z)f(y)|=2\pi.
\end{equation}
From (\ref{eqn3.2})-(\ref{eqn3.4}), we derive that
$|f(x)f(y)|=|xy|$. \hfill$\Box$

\vskip2mm

\noindent{\it Proof of the former part of Theorem D.}

We will give the proof by the induction on $k$. We first
note that $|q_iq_j|=\frac{\pi}{2}$ for any $i\neq j$ because
$\{q_1,\cdots,q_k\}$ is a $\frac\pi2$-separated subset and
$\diam(M)\leq \frac{\pi}{2}$. Hence, Theorem D is obviously true
when $k=1$ and 2. Now we assume that $k\geq3$.

By induction, there exists an isometrical embedding $g:
\Delta^{k-2}_+\to M$ such that $q_1,\cdots,q_{k-1}$ are the vertices
of $g(\Delta^{k-2}_+)$. We denote $g(\Delta^{k-2}_+) $ by $\Delta$,
and denote by $\Delta^\circ$ the interior part of $\Delta$. By
Theorem \ref{tct}, $\Delta$ belongs to $\{q_k\}^{\geq
\frac{\pi}{2}}$ (note that $q_1,\cdots,q_{k-1}$ are the vertices of
$\Delta$), which is convex in $M$ (Lemma \ref{lem2.2}). On the other
hand, $\{q_k\}^{\geq \frac{\pi}{2}}=\{q_k\}^{=\frac{\pi}{2}}$
because $\diam(M)\leq \frac{\pi}{2}$.  Then by Theorem \ref{tct=},
there is a triangle $\triangle q_kxy$ which is isometric to its
comparison triangle (in $\mathbb{S}^2$) for any $x,y\in \Delta$. (In
fact, the triangle $\triangle q_kxy$ bounds a convex domain which is
isometric to the convex domain bounded by the comparison triangle of
$\triangle q_kxy$ in $\mathbb{S}^2$ ([GM]); so if $k=3$, then the
proof is done.)

Now we fix a point $p\in \Delta^\circ$ and a geodesic $[q_kp]$. Note
that for any $x\in \Delta$, there is a unique geodesic $[px]$
between $p$ and $x$ in $M$ and $[px]\subset\Delta$ because $g$ is an
isometrical embedding. From the above we know that there is a
geodesic $[q_kx]$ such that the triangle $\triangle q_kpx$ composed
by $[q_kp], [q_kx]$ and $[px]$ is isometric to its comparison
triangle. Then we can define a map
$$h:\Delta\to \Sigma_{q_k}M \text{ by } x\mapsto
\uparrow_{q_k}^x.$$ Note that for any $x\in\Delta$ with $x\neq p$,
there is a unique $y\in\partial\Delta$ such that $x\in
[py]\subset\Delta$. Since $\triangle q_kpy$ bounds a convex domain
which is isometric to the convex domain bounded by the comparison
triangle of $\triangle q_kpy$ in $\Bbb S^2$, we can select $[q_kx]$
(for all $x\in\Delta$) such that $h([px])$ is a geodesic
$[\uparrow_{q_k}^p\uparrow_{q_k}^x]$ in $\Sigma_{q_k}M$. Hence, $h$
naturally induces a `tangential' map
$${\rm D}h:\Sigma_p\Delta\to \Sigma_{\uparrow_{q_k}^p}(\Sigma_{q_k}M)
\text{ defined by } \uparrow_p^x\mapsto
\uparrow_{\uparrow_{q_k}^p}^{\uparrow_{q_k}^x}.$$ Note that
$\Sigma_p\Delta=\mathbb{S}^{k-3}$ and
$\Sigma_{\uparrow_{q_k}^p}(\Sigma_{q_k}M)\in {\cal A}^{n-2}(1)$.

{\bf Claim}: {\it ${\rm D}h$ is a noncontractive map.} Since
triangles $\triangle q_kpx$ and $\triangle q_kpy$ are isometric to
their comparison triangles respectively for any $x,y\in \Delta$, we
have
\begin{equation}\label{eqn3.5}|\uparrow_{q_k}^p\uparrow_{q_k}^x|=|px| \text{ and }
\ |\uparrow_{q_k}^p\uparrow_{q_k}^y|=|py|.
\end{equation}
On the other hand, by Theorem \ref{tct} (on the triangle $\triangle
q_kxy$) we have
\begin{equation}\label{eqn3.6}
|\uparrow_{q_k}^x\uparrow_{q_k}^y|\geq|xy|
\end{equation}
(note that $|q_kx|=|q_ky|=\frac\pi2$). Then by the definition of
angles ([BGP]), we have
\begin{equation}\label{eqn3.7}
|\uparrow_{\uparrow_{q_k}^p}^{\uparrow_{q_k}^x}\uparrow_{\uparrow_{q_k}^p}^{\uparrow_{q_k}^y}|
\geq |\uparrow_p^x\uparrow_p^y|,
\end{equation}
i.e. the claim is verified.

By Lemma \ref{lem3.1}, the claim implies that ${\rm D}h$ is an
isometrical embedding in fact, so the inequality (\ref{eqn3.7}) is
an equality:
\begin{equation}\label{eqn3.8}
|\uparrow_{\uparrow_{q_k}^p}^{\uparrow_{q_k}^x}\uparrow_{\uparrow_{q_k}^p}^{\uparrow_{q_k}^y}|
=|\uparrow_p^x\uparrow_p^y|.
\end{equation}
Note that (\ref{eqn3.5}) and (\ref{eqn3.8}) imply that the hinge
$p\prec_x^y\ \subset \Delta\subset\mathbb{S}^{k-2}$ is the
comparison hinge of the hinge
$\uparrow_{q_k}^p\prec_{\uparrow_{q_k}^x}^{\uparrow_{q_k}^y}\subset
\Sigma_{q_k}M$. Then by Theorem \ref{tct}, we have
$|\uparrow_{q_k}^x\uparrow_{q_k}^y|\leq|xy|$, which together with
(\ref{eqn3.6}) implies that
\begin{equation}\label{eqn3.9}
|\uparrow_{q_k}^x\uparrow_{q_k}^y|=|xy|.
\end{equation}
Note that there is a unique geodesic $[xy]$ between any
$x\in\Delta^\circ$ and $y\in\Delta^\circ$. By Theorem \ref{tct=},
(\ref{eqn3.9}) implies that the triangle $\triangle q_kxy$ composed
by geodesics $[q_kx]$, $[q_ky]$ and $[xy]$ is isometric to its
comparison triangle (in $\mathbb{S}^{2}$). Therefore, we can
conclude that $\{q_k\}*\Delta^\circ$ can be isometrically embedded
into $M$, so $\{q_k\}*\Delta$ can be isometrically embedded into $M$
(by a standard limit argument). Recall that
$\Delta=g(\Delta^{k-2}_+)$, whose vertices are $q_1,\cdots,q_{k-1}$.
It then follows that $\Delta^{k-1}_+$ $(=\{q_k\}*\Delta)$ can be
isometrically embedded into $M$ with $q_1,\cdots,q_{k}$ being the
vertices. \hfill $\Box$

\vskip2mm

\noindent{\it Proof of the latter part of Theorem D under the
assumption $k=n+1$.}

We will give the proof by the induction on $n$.

Obviously, when $n=1$, $M$ is either an arc of length $\frac\pi2$
with $q_1$ and $q_2$ being end points or a circle of perimeter $\pi$
with $q_1$ and $q_2$ being antipodal points.

Now we assume that $n>1$. From the proof of the former part of
Theorem D, we know that
$N\triangleq\{q_{n+1}\}^{\geq\frac\pi2}=\{q_{n+1}\}^{=\frac\pi2}$
which is convex in $M$; and thus $N\in{\cal A}(1)$ with
$\diam(N)\leq\frac\pi2$, and $\dim(N)\leq n-1$ (Lemma \ref{lem2.5}).
On the other hand, note that $\{q_1,\cdots,q_n\}\subset N$, so by
the former part of Theorem D we have $\dim(N)=n-1$. Hence, by
induction we can conclude that $N$ is a glued space of finite copies
of $\Delta^{n-1}_+$ along some ``faces'' $\Delta^{n-2}_+$ of them.

Given an arbitrary $\Delta^{n-1}_+\subset N$ and any point $p\in
(\Delta^{n-1}_+)^\circ$, from the above proof for the former part,
we know that a geodesic $[q_{n+1}p]$ determines a $\Delta^{n}_+$
(with $[q_{n+1}p]\subset\Delta^{n}_+$ and $q_1,\cdots,q_{n+1}$ being
the vertices) which can be isometrically embedded into $M$. As a
result, $\Sigma_p\Delta^{n}_+$ can be isometrically embedded into
$\Sigma_pM$ ([BGP]). Note that
$$\Sigma_p\Delta^{n}_+=\{\uparrow_p^{q_{n+1}}\}*\Sigma_p\Delta^{n-1}_+=
\{\uparrow_p^{q_{n+1}}\}*\Bbb S^{n-2}$$ (which is a half $\Bbb
S^{n-1}$). This implies that  if there is another geodesic
$[q_{n+1}p]'$, then
$\Sigma_pM=\{\uparrow_p^{q_{n+1}},(\uparrow_p^{q_{n+1}})'\}*\Bbb
S^{n-2}=\Bbb S^{n-1}$. Hence, there are at most two geodesics
between $p$ and $q_{n+1}$, and so there are at most two
$\Delta^{n}_+$ (with $q_1,\cdots,q_{n+1}$ being the vertices) which
contain the given $\Delta^{n-1}_+$.  On the other hand, of course,
every $\Delta^{n}_+\subset M$ with $q_1,\cdots,q_{n+1}$ being the
vertices contains a $\Delta^{n-1}_+$ in $N$. It follows that $M$
contains only finite copies of $\Delta^{n}_+$ because $N$ is a glued
space of finite copies of $\Delta^{n-1}_+$.

Now we let $M'$ denote the union of all $\Delta^{n}_+$ with
$q_1,\cdots,q_{n+1}$ being the vertices.  {\bf Claim}: $M'=M$. If
the claim is not true, then for any $x\in M\setminus M'$ there is
$p\in M'$ such that $|xp|=\min_{q\in M'}\{|xq|\}$ (note that $M'$ is
compact). Obviously, $p$ cannot be an interior point of any
$\Delta^{n}_+$ in $M'$.

\noindent{\bf Subclaim}: {\it $p$ cannot be an interior point of any
``face'' of any $\Delta^{n}_+$ in $M'$ either}. If the subclaim is
not true, then we can rearrange all $q_i$ such that $p$ is an
interior point of $\Delta^{n-1}_+\triangleq\Delta^{n}_+\cap N$ for
some $\Delta^{n}_+$ in $M'$.  Let $[q_{n+1}p]$ be the geodesic
between $q_{n+1}$ and $p$ in the $\Delta^{n}_{+}$. Note that
$\Sigma_p\Delta^{n-1}_{+}$ ($=\Bbb S^{n-2}$) is  convex in
$\Sigma_pM$, and $\Sigma_p\Delta^{n}_{+}=\{\uparrow_p^{q_{n+1}}\}$
$*\Sigma_p\Delta^{n-1}_{+}$ (which is a half $\Bbb S^{n-1}$). By
Lemma \ref{lem2.3}, we have
$(\Sigma_p\Delta^{n-1}_{+})^{\geq\frac\pi2}=(\Sigma_p\Delta^{n-1}_{+})^{=\frac\pi2}$.
Now we select a geodesic $[xp]$. Since
$|xp|=\min_{q\in\Delta^{n}_{+}}\{|xq|\}$, by the first variation
formula ([BGP]) we have $|\uparrow_p^x\xi|\geq\frac\pi2$ for any
$\xi\in\Sigma_p\Delta^{n-1}_{+}$, i.e. $\uparrow_p^x\in
(\Sigma_p\Delta^{n-1}_{+})^{\geq\frac\pi2}$. It then follows that
$\{\uparrow_p^x,\uparrow_p^{q_{n+1}}\}\subset
(\Sigma_p\Delta^{n-1}_{+})^{=\frac\pi2}$. On the other hand,
$(\Sigma_p\Delta^{n-1}_{+})^{\geq\frac\pi2}$ is convex in
$\Sigma_pM$ (Lemma \ref{lem2.2}); then
$\dim((\Sigma_p\Delta^{n-1}_{+})^{=\frac\pi2})=0$ (Lemma
\ref{lem2.5}), and so
$(\Sigma_p\Delta^{n-1}_{+})^{=\frac\pi2}=\{\uparrow_p^x,\uparrow_p^{q_{n+1}}\}$
with $|\uparrow_p^x\uparrow_p^{q_{n+1}}|=\pi$. By the first
variation formula ([BGP]), there exists another geodesic
$[pq_{n+1}]'$ in $M$ such that
$|\uparrow_p^x(\uparrow_p^{q_{n+1}})'|\leq\frac\pi2$ (otherwise
there is $x'\in[px]$ such that $|q_{n+1}x'|>\frac\pi2$ which
contradicts $\diam(M)\leq\frac\pi2$). However, from the above,
$[pq_{n+1}]'$ determines another $\Delta^{'n}_+$ with
$q_1,\cdots,q_{n+1}$ being the vertices and
$[pq_{n+1}]'\subset\Delta^{'n}_+$, and thus
$(\uparrow_p^{q_{n+1}})'\in(\Sigma_p\Delta^{n-1}_{+})^{=\frac\pi2}$.
Hence, $(\uparrow_p^{q_{n+1}})'$ has to be $\uparrow_p^{x}$, and so
$[px]\subset[pq_{n+1}]'\subset\Delta^{'n}_+$ which contradicts $x\in
M\setminus M'$ (i.e., the subclaim is verified).

For convenience, a ``face'' $\Delta^{n-1}_+$ of any
$\Delta^{n}_+\subseteq M'$ (resp. a ``face'' $\Delta^{n-2}_+$ of
such a $\Delta^{n-1}_+$) is said to be an $(n-1)$-face (resp. $(n-2)$-face)
of $M'$. According to the subclaim and its proof, it is not hard to
observe that:

\vskip1mm

\noindent(3.10)\ {\it If $M\neq M'$, then for any $x\in M\setminus
M'$ and any interior point $q$ of any $\Delta^{n}_+$ in $M'$, the
nearest point in $ [xq]\cap M'$ to $x$ for any $[xq]$ has to lie in
an $(n-2)$-face of $M'$.}

\vskip1mm

\noindent However, by the induction on $n$, we will derive a
contradiction under (3.10). If $n=2$, note that (3.10) implies that
there must be a geodesic which branches at some $0$-face (which is a
point) of $M'$, which is impossible ([BGP]). Now we assume that
$n>2$. Let $v$ be the nearest point to $x$ on $[xq]\cap M'$ (where
$[xq]$ is the geodesic in (3.10)). Note that $\Sigma_vM'$ is also the
union of finite copies of $\Delta^{n-1}_+$ and each $\Delta^{n-1}_+$
can be isometrically embedded into $\Sigma_vM$ (because each
$\Delta^{n}_+$ of $M'$ can be isometrically embedded into $M$); and
note that $\Sigma_vM\neq\Sigma_vM'$ because
$\uparrow_v^x\not\in\Sigma_vM'$ (note that $[vx]\cap M'=\{v\}$).
Then from the definition of the angle (i.e. the distance between two
directions in $\Sigma_vM$) ([BGP]), it is not hard to see that (3.10)
implies that:

\vskip1mm

\noindent(3.11)\ {\it  For any $\xi\in\Sigma_vM\setminus\Sigma_vM'$
and any interior point $\eta$ of any $\Delta^{n-1}_+$ in
$\Sigma_vM'$, the nearest point in $[\xi\eta]\cap\Sigma_vM'$ to $\xi$ for
any $[\xi\eta]$ has to lie in an $(n-3)$-face of $\Sigma_vM'$.}

\vskip1mm

\noindent By induction, we can derive a contradiction under (3.11).
That is, we get a contradiction under (3.10), so we have $M=M'$ (i.e.
the claim is verified).

Now we can conclude that $M$ is not only the union of all
$\Delta^{n}_+$ with $q_1,\cdots,q_{n+1}$ being the vertices by the
claim, but also a glued space of these $\Delta^{n}_+$ along some
$(n-1)$-faces (not $(n-k)$-faces with $k\geq2$) of them from the
proof of the claim.\hfill $\Box$

\vskip2mm

\setcounter{equation}{11}

\begin{remark}\label{rem3.12} \rm In fact, if $k=n+1$ in Theorem D, then
we can determine all the possible structures of $M$ by the induction
on the dimension $n$. From the above proof, we know that
$N\triangleq\{q_{n+1}\}^{=\frac\pi2}$ belongs to ${\cal A}^{n-1}(1)$
(because $N$ is convex in $M$) with $\diam(N)\leq\frac\pi2$ and
$\{q_1,\cdots,q_n\}\subset N$; so by induction we can determine all
the possible structures of $N$. Moreover, we know that there are at
most two geodesics between $q_{n+1}$ and any interior point of any
$\Delta^{n-1}_+$ in $N$. On the other hand, $M$ is a glued space of
all $\Delta^{n}_+$ along some $(n-1)$-faces of them. Hence, there
are at most two geodesics between $q_{n+1}$ and any point in $N$;
and for any $x\in M$, there is some geodesic $[q_{n+1}p]$ with $p\in
N$ such that $x\in[q_{n+1}p]$. Then we can prove that either
$M=\{q_{n+1}\}*N$, or $\tilde N\triangleq\Sigma_{q_{n+1}}M$ admits
an isometry $\sigma$ of order 2 (i.e. $\sigma^2=\text{id}$), which
naturally induces an isometry $\tilde\sigma$ of order 2 on the
suspension $\{q_{n+1},\bar q_{n+1}\}*\tilde N$ (where $|q_{n+1}\bar
q_{n+1}|=\pi$) with $\tilde\sigma([q_{n+1}p])=[\bar
q_{n+1}\sigma(p)]$ and $\tilde\sigma([\bar q_{n+1}p])=[q_{n+1}\sigma
(p)]$, such that $N=\tilde N/\langle\sigma\rangle$ and
$M=(\{q_{n+1},\bar q_{n+1}\}*\tilde N)/\langle\tilde\sigma\rangle$
(one can give the detailed proof for this by referring to the proof
of Lemma \ref{lem2.5} in [SSW], or [RW]). This implies that we can
determine the structure of $M$ by that of $N$. For  example, we give
all the possible structures of $M$ for $n=1$, $2$ and $k=n+1$ (for
convenience, we let $\Bbb Z_2$ denote both $\langle\sigma\rangle$
and $\langle\tilde\sigma\rangle$).

\vskip1mm

\noindent $n=1$ and $k=2$: $M=\Delta^{1}_+$ $(=[q_1q_2])$, an arc of
length $\frac{\pi}{2}$; or $M=S^1_{\pi}$,  a circle of perimeter
$\pi$, which is a glued space of two copies of $\Delta^{1}_+$  at
$q_1$ and $q_2$.

\vskip1mm

\noindent $n=2$ and $k=3$: $M=\Delta^{2}_+$ $(=\{q_3\}*[q_1q_2])$,
$\{q_3\}*S^1_{\pi}$ (a glued space of two $\Delta^{2}_+$ along
$[q_3q_1]$ and $[q_3q_2]$), $(\{q_{3},\bar q_{3}\}*S^1_{\pi})/\Bbb
Z_2$ (a glued space of two $\Delta^{2}_+$ along $[q_3q_1]$,
$[q_3q_2]$ and $[q_1q_2]$) where $\Bbb Z_2$ acts on $S^1_{\pi}$ by a
reflection (note that $S^1_{\pi}/\Bbb Z_2=[q_1q_2]$), or
$(\{q_{3},\bar q_{3}\}*S^1_{2\pi})/\Bbb Z_2$ (=${\Bbb R\Bbb P}^2$,
a glued space of four $\Delta^{2}_+$ along their boundaries) where
$\Bbb Z_2$ acts on $S^1_{2\pi}$ by the antipodal map (note that
$S^1_{2\pi}/\Bbb Z_2=S^1_{\pi}$).

\end{remark}

\vskip2mm

We will end this section by giving a brief proof for Corollary E.

\vskip1mm

\noindent{\bf Proof of Corollary E}:

Note that $M$ has  empty boundary (because $M$ is closed). According
to Remark \ref{rem3.12}, it has to hold that $\tilde N$
$(=\Sigma_{q_{n+1}}M)$ admits an isometrical $\Bbb Z_2$-action which
naturally induces an isometrical $\Bbb Z_2$-action on the suspension
$\{q_{n+1},\bar q_{n+1}\}*\tilde N$ such that $M=(\{q_{n+1},\bar
q_{n+1}\}*\tilde N)/\Bbb Z_2$. Since $M$ is a Riemannian manifold,
we have $\tilde N=\Sigma_{q_{n+1}}M=\Bbb S^{n-1}$; and thus
$\{q_{n+1},\bar q_{n+1}\}*\tilde N=\Bbb S^{n}$, so $M=\Bbb
S^{n}/\Bbb Z_2$. Since $\Bbb Z_2$ acts on $\Bbb S^n$ by isometries
and $M$ is a Riemannian manifold (with $\diam(M)\leq\frac\pi2$), it
has to hold that $M={\Bbb R\Bbb P}^n$ (i.e. the $\Bbb Z_2$-action on
$\Bbb S^{n}$ must be realized by the antipodal map). \hfill$\Box$

%%%%%%%%%%%%%%%%%%%%%%%%%%%%%% % Section 4 %%%%%%%%%%%%%%%%%%%%%%%%%%%%%%%%%%%%%%%

\section{A technical corollary of Theorem D}

We first give an easy corollary of Theorem D.

\begin{prop}\label{prop4.1}
Let $M\in {\cal A}^n(1)$, and let $\{p_1,\cdots,p_k\}$ be a
$\frac\pi2$-separated subset in $M$ with $|p_1p_i|>\frac{\pi}{2}$ for
$2\leq i\leq k$. Suppose that $N$ is a complete and convex subset in
$\{p_1,\cdots,p_k\}^{\geq\frac\pi2}$ with $\diam(N)\leq
\frac{\pi}{2}$, and that $\{q_1,\cdots,q_h\}$ is a
$\frac\pi2$-separated subset in $N$. Then $k+h\leq n+2$.
\end{prop}

\noindent{\it Proof.}  Note that $N\in {\cal A}(1)$ (because $N$ is
convex in $\{p_1,\cdots,p_k\}^{\geq\frac\pi2}$ which is convex in
$M$ by Lemma \ref{lem2.2}) with $\diam(N)\le\frac\pi2$. By Theorem
D, $\Delta^{h-1}_+$ can be isometrically embedded into $N$ with
$q_1,\cdots,q_{h}$ being the vertices. Now we consider
$\Sigma_{q_h}M\in {\cal A}^{n-1}(1)$. Similar to the proof of Theorem A, we can conclude that any
$\{\uparrow_{q_h}^{p_1},\cdots,\uparrow_{q_h}^{p_k}\}$ is a
$\frac\pi2$-separated subset in $\Sigma_{q_h}M$ with
$|\uparrow_{q_h}^{p_1}\uparrow_{q_h}^{p_i}|>\frac{\pi}{2}$ for
$2\leq i\leq k$. Moreover,
$\Sigma_{q_h}\Delta^{h-1}_+=\Delta^{h-2}_+$, which can be
isometrically embedded into
$\{\uparrow_{q_h}^{p_1},\cdots,\uparrow_{q_h}^{p_k}\}^{\geq\frac\pi2}$
(with some $\uparrow_{q_h}^{q_1},\cdots,\uparrow_{q_h}^{q_{h-1}}$
being the vertices). Then we can repeat such an argument on
$\Sigma_{\uparrow_{q_h}^{q_{h-1}}}(\Sigma_{q_h}M)$, so by induction
(note that when $h=0$, the proposition is obvious by Theorem
A) we obtain that $k+h-1\leq n-1+2$, i.e. $k+h\leq n+2$. \hfill
$\Box$

\vskip2mm

Based on Proposition \ref{prop4.1}, we give a technical corollary
which may make sense in analyzing the direction spaces as same as in
[P1].

\begin{coro}\label{cor4.2} {\rm \bf (A technical corollary of Theorem D)}
Let $M$, $\{p_1,\cdots,p_k\}$ and $N$ be the same as in Proposition
\ref{prop4.1}. Then given any $z_1\in N$ and sufficiently small
$\epsilon>0$, there exist $\{z_2,\cdots,z_h\}\subset\partial
B(z_1,\epsilon)\cap N$ with $h\le n+2-k$ and $|z_iz_j|>\epsilon$
($2\le i\neq j\le h$) and $\delta_1, \delta_2>0$ with
$h\delta_2<\delta_1$ such that for any $p\in M$ either
$|pN|>\delta_2$ or there is some $1\le i\le h$ such that
$|pz_i|<\frac{\pi}{2}-\delta_1$ and
$|pz_{i'}|<\frac{\pi}{2}+\delta_2$ for all $i'\ne i$.
\end{coro}

\noindent{\it Proof.} Note that we can find a maximal
$\frac\pi2$-separated subset in $N$, $\{z_1, q_2,\cdots,q_h\}$ (i.e.
$|x\{z_1, q_2,\cdots,q_h\}|<\frac\pi2$ for all $x\in N$). By
Proposition \ref{prop4.1}, we have $h\le n+2-k$. In $N$, we consider
$L\triangleq\{z_1\}^{\geq\frac\pi2}$, which is convex in $N$ (so in
$M$) by Lemma \ref{lem2.2} (so $L\in {\cal A}(1)$). Since
$\diam(N)\leq\frac\pi2$, we have $L=\{z_1\}^{=\frac\pi2}$ and
$\{q_2,\cdots,q_h\}\subset L$.

Next we will find $z_i$ for $2\le i\le h$. Note that we can select
$\bar q_i\in L^\circ$ (the interior part of $L$) such that $|\bar
q_iq_i|\ll\epsilon$ ([BGP]) (if $q_i\in L^\circ$, then we let $\bar
q_i$ just be $q_i$). Then we select an arbitrary geodesic $[z_1\bar
q_i]$, and select $z_i\in [z_1\bar q_i]$ with $|z_1z_i|=\epsilon$.
Of course, $\{z_2,\cdots,z_h\}\subset\partial B(z_1,\epsilon)\cap
N$. By Theorem \ref{tct} on the triangle $\triangle z_1\bar q_i\bar
q_j$, it is not hard to see that $|z_iz_j|>\epsilon$ for $2\le i\neq
j\le h$ (note that $|z_1\bar q_i|=|z_1\bar q_j|=\frac\pi2$ and
$|\bar q_i\bar q_j|$ is almost equal to $\frac\pi2$, and $z_i\in
[z_1\bar q_i],z_j\in [z_1\bar q_j]$).

\noindent{\bf Claim}: {\it There exists $\delta>0$ such that
$|x\{z_1,\cdots,z_h\}|<\frac\pi2-\delta$ for any $x\in N$ (i.e.
there is some $1\le i\le h$ such that $|xz_i|<\frac\pi2-\delta$).}

It is easy to see that the claim implies that we can find the
desired $\delta_1$ and $\delta_2$, so we only need to verify the
claim in the rest of proof.

Since $\{z_1, q_2,\cdots,q_k\}$ is a maximal $\frac\pi2$-separated
subset in $N$ and $|\bar q_iq_i|\ll\epsilon$, there is
$\epsilon_1>0$ such that $|y\{\bar q_2,\cdots,\bar
q_k\}|<\frac\pi2-\epsilon_1$ for any $y\in L$, i.e. there is a $\bar
q_i$ such that $|y\bar q_i|<\frac\pi2-\epsilon_1$. On the other
hand, by the first variation formula ([BGP]) we have
$|\uparrow_{\bar q_i}^{z_1}\xi|\geq\frac\pi2$ for any
$\xi\in\Sigma_{\bar q_i}L$ (note that $L=\{z_1\}^{=\frac\pi2}$ which
is convex in $N$). In fact, by Lemma \ref{lem2.3} we have
$|\uparrow_{\bar q_i}^{z_1}\xi|=\frac\pi2$ because $\Sigma_{\bar
q_i}L$ is convex in $\Sigma_{\bar q_i}N$ and $\Sigma_{\bar q_i}L$
has empty boundary (note that $L$ is convex in $N$ and $\bar q_i\in
L^\circ$). It then follows that $|\uparrow_{\bar
q_i}^{z_i}\uparrow_{\bar q_i}^{y}|=\frac\pi2$ for any geodesic
$[\bar q_iy]\subset L$. And thus by Theorem \ref{tct} on the hinge
$\bar q_i\prec^{z_i}_y$, there is $\chi(\epsilon,\epsilon_1)>0$
(where $\chi(\epsilon,\epsilon_1)\to0$ as $\epsilon,\epsilon_1\to
0$) such that
$$|yz_i|<\frac\pi2-\chi(\epsilon,\epsilon_1)$$
(note that $|\bar q_iy|<\frac\pi2-\epsilon_1$, $|\bar
q_iz_i|=\frac\pi2-\epsilon$ and $|\uparrow_{\bar
q_i}^{z_i}\uparrow_{\bar q_i}^{y}|=\frac\pi2$). Hence, for any $x\in
B(L,\epsilon_2)\cap N$ (the $\epsilon_2$-tubular neighborhood of $L$
in $N$) there exist a $z_i$ and $\chi(\epsilon,\epsilon_1,\epsilon_2)>0$ such that
$$|xz_i|<\frac\pi2-\chi(\epsilon,\epsilon_1,\epsilon_2).$$
On the other hand, there is a $\epsilon_3>0$ such that
$|xz_1|<\frac\pi2-\epsilon_3$ for any $x\in N\setminus
B(L,\epsilon_2)$ (otherwise there is $x_0\in N\setminus
B(L,\epsilon_2)$ such that $|x_0z_1|=\frac\pi2$, i.e. $x_0\in L$; a
contradiction). Note that
$\delta\triangleq\min\{\chi(\epsilon,\epsilon_1,\epsilon_2),\epsilon_3\}$
is the desired number of the claim. \hfill $\Box$

%%%%%%%%%%%%%%%%%%%%%%%%%%%%%%%%%%%%%%%%%%%%%%%%%%%%%%%%%%%%%%%%%%%%%

\vskip1cm

\noindent{\bf \Large Appendix}

\vskip5mm

\noindent{\bf Proof of Theorem \ref{thm0.2}:}

We will give the proof by the induction on the dimension $n$.
Obviously, Theorem \ref{thm0.2} is true if $n=0$ and 1. Now we
assume $n>1$. We consider $\Sigma_{q_1}M$ which belongs to ${\cal
A}^{n-1}(1)$ ([BGP]). By Lemma \ref{lem1.2}, it is easy to see that
in $\{q_2,\cdots,q_k\}$ there is at most one point, say $q_k$, such
that $|q_1q_k|=\pi$. Then by Lemma \ref{lem1.1}, any
$\{\uparrow_{q_1}^{q_2},\cdots,\uparrow_{q_1}^{q_{k-1}}\}$ is a
$\frac\pi2$-separated subset in $\Sigma_{q_1}M$. By the inductive
assumption on $\Sigma_{q_1}M$, we have $$k-2\leq2(n-1+1), \text{
i.e., } k\leq2(n+1).$$ Moreover, if $k=2(n+1)$, then it has to hold
that $|q_1q_{2(n+1)}|=\pi$; and so $M=\{q_1,q_{2(n+1)}\}*M_1$ for
some $M_1\in {\cal A}^{n-1}(1)$ (Lemma \ref{lem1.2}). This implies
that $\{q_2,\cdots,q_{2n+1}\}$ is a $\frac\pi2$-separated subset in
$M_1$. By the inductive assumption on $M_1$, we have that $M_1$ is
isometric  to the unit sphere $\mathbb{S}^{n-1}$. Therefore, we can
conclude that $k=2(n+1)$ if and only if $M$ is isometric  to
$\{q_1,q_{2(n+1)}\}*\mathbb{S}^{n-1}$ which is the unit sphere
$\mathbb{S}^{n}$. \hfill$\Box$

\vskip2mm

\noindent{\bf Proof of Lemma \ref{lem1.4}:}

We will give the proof by the induction on $n$. Obviously, if $n=1$,
the lemma is true (because $M$ is an arc of length $\leq \pi$ if
$M\in {\cal A}^1(1)$ has nonempty boundary (see the convention after
Lemma \ref{lem1.2})). Now we assume that $n>1$.

Since $\partial M$ is compact (because $\partial M$ is closed in $M$
([BGP]) and $M$ is compact), we select $q\in\partial M$ with
$|pq|=|p\partial M|$ and a geodesic $[pq]$. By the first variation
formula,
$$|\uparrow_q^p\xi|\geq\frac\pi2$$
in $\Sigma_qM\in{\cal A}^{n-1}(1)$ for any
$\xi\in\partial(\Sigma_qM)$ (refer to [BGP] for the details on the
boundary in Alexandrov geometry). By induction, we have
$\Sigma_qM=\{\uparrow_q^p\}*\partial(\Sigma_qM)$, so
$|\uparrow_q^p\eta|\leq\frac\pi2$ for any $\eta\in\Sigma_qM$ and the
``='' holds if and only if $\eta\in\partial(\Sigma_qM)$. Then by
Theorem \ref{tct} on any triangle $\triangle pqr$ with $r\in\partial
M$ (note that $|pr|\geq|pq|\ge\frac\pi2$ and $\angle
pqr\leq\frac\pi2$), it has to hold that
$$|pr|=|pq|=\frac\pi2 \text{ and } \angle pqr=\frac\pi2.$$
It then follows that $|pr|=|p\partial M|=\frac\pi2$. And
$\uparrow_q^r\in\partial(\Sigma_qM)$ for any geodesic $[qr]$, so
$[qr]$ belongs to $\partial M$ ([BGP]). Now we can take the place of
$q$ by any $r\in\partial M$. Similarly, we can get that any geodesic
$[rr']$ for any other $r'\in\partial M$ belongs to $\partial M$,
i.e. $\partial M$ is convex in $M$; and
$\Sigma_rM=\{\uparrow_r^p\}*\partial(\Sigma_rM)$, which implies that
there is a unique geodesic between $p$ and $r$. Hence, in order to
prove that $M=\{p\}*\partial M$, it suffices to show that there is a
point $r$ in $\partial M$ such that $x\in [pr]$ for any $x\in M$
(cf. Proposition \ref{prop2.8} and refer to [SSW]). Similarly, we
can select $s\in\partial M$ such that $|xs|=|x\partial M|$ and a
geodesic $[xs]$; and get that $|\uparrow_s^x\xi|\geq\frac\pi2$ in
$\Sigma_sM$ for any $\xi\in\partial(\Sigma_sM)$, and so
$\Sigma_sM=\{\uparrow_s^x\}*\partial(\Sigma_sM)$. It then has to
hold that $$\uparrow_s^x=\uparrow_s^p.$$ On the other hand, note
that we have proved that $|x\partial M|\leq\frac\pi2$. It then
follows that $[xs]\subseteq[ps]$ (of course $x\in [ps]$).
\hfill$\Box$

\vskip2mm

\noindent{\bf Proof of Lemma \ref{lem2.11}:}

Due to Proposition \ref{prop2.8}, it suffices to show that there is a unique geodesic between any $x\in X$ and $y\in S^1$.
Note that $\Sigma_xX$  is a complete convex subset in $\Sigma_xM$ and $\Sigma_xX$ has empty boundary because
$X$ is a complete convex subset in $M$ and $X$ has empty boundary. By Lemma \ref{lem2.2}, $(\Sigma_xX)^{\geq\frac\pi2}$ is
convex in $\Sigma_xM$; and by Lemma \ref{lem2.3}, we have $(\Sigma_xX)^{\geq\frac\pi2}=(\Sigma_xX)^{=\frac\pi2}$. It then follows
from Lemma \ref{lem2.5} that $\dim((\Sigma_xX)^{=\frac\pi2})\leq 1$ (note that $\dim(\Sigma_xX)=n-3$ and $\dim(\Sigma_xM)=n-1$).
Moreover, for any $y\in S^1$ and any geodesic $[xy]$, by the first variation formula we have that $\uparrow_x^y$
belongs to $(\Sigma_xX)^{\geq\frac\pi2}$ ($=(\Sigma_xX)^{=\frac\pi2}$). This implies that $\dim((\Sigma_xX)^{=\frac\pi2})=1$.
On the other hand, since $S^1\subset X^{=\frac\pi2}$ and $S^1$ is convex in $M$,
by Theorem \ref{tct=}, for any geodesic $[y'y'']\subset S^1$ with $y\in[y'y'']$ there is a triangle $\triangle xy'y''$ containing $[y'y'']$ which
is isometric to its comparison triangle; moreover, $[xy]$ belongs to the convex domain bounded by
$\triangle xy'y''$ ([GM]). Then we can define a multi-value map
$$f:S^1\to (\Sigma_xX)^{=\frac\pi2} \text{ by } y\mapsto \{\text{all directions from $x$ to $y$}\}$$
such that for any $y\in S^1$ and any geodesic $[xy]$ there are a neighborhood $U$ of $y$ in $S^1$ and
a neighborhood $\tilde U$ of $\uparrow_x^y$ in $(\Sigma_xX)^{=\frac\pi2}$ such that $f|_U:U\to \tilde U$ is an isometry.
This implies that $f(S^1)$ is both open and closed in $(\Sigma_xX)^{=\frac\pi2}$, so we have $f(S^1)=(\Sigma_xX)^{=\frac\pi2}$.
It therefore follows that $f^{-1}$ is a covering map. Since $(\Sigma_xX)^{=\frac\pi2}\in {\cal A}^1(1)$ and $S^1$ has perimeter
$>\pi$, it has to hold that $f^{-1}$ is a 1-1 covering map (of course  $f$ is also a 1-1 map),
which implies that there is a unique geodesic between $x\in X$ and $y\in S^1$.
\hfill$\Box$

\vskip2mm

\noindent{\bf Remark A.1.}\ \  According to the proof of Lemma \ref{lem2.11}, if $X$ has nonempty boundary in Lemma \ref{lem2.11},
we can still conclude that there is a unique geodesic between any interior point $x$ in $X$ and $y\in S^1$; as a result, $X*S^1$
can be isometrically embedded into $M$ (cf. Proposition \ref{prop2.8} and refer to [SSW]).

\vskip2mm

\noindent{\bf Remark A.2.}\ \ In fact, Lemma \ref{lem2.11} has the following
generalized version: {\it Let $X$ and $Y$ be two complete convex subsets in $M\in{\cal
A}^n(1)$. Then $M=X*Y$ if\ \ {\rm (i)} $\dim(X)+\dim(Y)=n-1$; {\rm (ii)} both
$X$ and $Y$ have empty boundary; {\rm (iii)} either $X$ or $Y$ has radius
$>\frac\pi2$; {\rm (iv)} $|xy|=\frac\pi2$ for any
$x\in X$ and $y\in Y$}.

\vskip1cm

\noindent{\bf Acknowledgements.}\ \ 
The authors would like to thank Professor Frederick Wilhelm for supplying the reference [GW].

%%%%%%%%%%%%%%%%%%%%%%%%%%%%%%%%%%%%%%%%%%%%%%%%%%%%%%%%%%%%%%%

\noindent School of Mathematical Sciences (and Lab. math. Com.
Sys.), Beijing Normal University, Beijing, 100875
P.R.C.\\
e-mail: suxiaole$@$bnu.edu.cn; wwyusheng$@$gmail.com

\vskip2mm

\noindent Mathematics Department, Capital Normal University,
Beijing, 100037 P.R.C.\\
e-mail: hwsun$@$bnu.edu.cn

\end{document}